\documentclass[12pt,a4paper]{article}
\usepackage{color}
\usepackage{amsmath}
\usepackage{amsfonts}
\usepackage{amssymb}
\usepackage{latexsym}
\usepackage{latexsym}
\usepackage{array}
\usepackage{hyperref}
\usepackage{amsthm}
\usepackage{graphics}

\newtheorem{Theorem}{Theorem}
\newtheorem{Lemma}{Lemma}

\voffset=-.5in \hoffset=-.5in \numberwithin{equation}{section}

\begin{document}
\title {
Multi-scaling Limits for Relativistic Diffusion Equations with Random Initial
Data
}


\author { Gi-Ren Liu\thanks{PhD Class in Mathematics, National Taiwan University, Taipei 10617,
Taiwan.} \quad \quad Narn-Rueih Shieh\thanks{ Mathematics
Department, Honorary Faculty, National Taiwan University, Taipei
10617, Taiwan.
 Correspondence to
E-mail: \texttt{shiehnr@ntu.edu.tw } URL:
\texttt{http://www.math.ntu.edu.tw/\~{}shiehnr/}} }

\date{}

 \maketitle
\vskip 18 pt





\begin{abstract}
{\small  Let $u(t,\mathbf{x}),\ t>0,\ \mathbf{x}\in
\mathbb{R}^{n},$ be the spatial-temporal random field arising from
the solution of a relativistic diffusion equation  with the
spatial-fractional  parameter $\alpha\in (0,2)$ and the mass
parameter $\mathfrak{m}> 0$, subject to a random initial condition
$u(0,\mathbf{x})$ which is characterized as a subordinated
Gaussian field. In this article, we study the large-scale  and the
small-scale  limits for the suitable space-time re-scalings   of
the solution field $u(t,\mathbf{x})$. Both the Gaussian and the
non-Gaussian limit theorems are discussed. The small-scale scaling
involves not only to scale on $u(t,\mathbf{x})$ but also to
re-scale the initial data; this is a new-type result for the
literature. Moreover, in the two scalings the parameter $\alpha\in
(0,2)$ and the parameter $\mathfrak{m}> 0$ paly distinct roles for
the scaling  and the limiting procedures. }
\end{abstract}


{\it Short title:} Multi-Scaling Limits for Relativistic Diffusion
Equations

 {\it 2010 AMS classification numbers:} 60G60; 60H05; 62M15;
      35K15.

 {\it Key words:} Large-scale limits; Small-scale limits; Relativistic diffusion equations;
 Random initial data; Multiple It$\hat{\textup{o}}$-Wiener integrals;
 Subordinated Gaussian fields; Hermite ranks.


\newpage

\section{Introduction}

In this paper, we consider  the following Cauchy problem for the
relativistic diffusion equation (RDE for brevity), {\em subject to
some random initial data}, and aim to discuss the scaling limits
for the spatial-temporal random field arising from the solution of
this random initial value problem :
\begin{align}\label{relativistic diffusion equation}
\frac{\partial}{\partial
t}u(t,\mathbf{x})=(\mathfrak{m}-(\mathfrak{m}^{\frac{2}{\alpha}}-\Delta)^{\frac{\alpha}{2}})
u(t,\mathbf{x}),
\ \ u(0,\mathbf{x})=u_{0}(\mathbf{x}),\ t\geq0,\ \mathbf{x}\in \mathbb{R}^{n},
\end{align}
 with the spatial-fractional  parameter $\alpha\in (0,2)$ and the  (normalized) mass parameter $\mathfrak{m}> 0$.

RDEs appear in vast literature of mathematics and physics. The
prominent case is $\alpha=1$, for which
$-(\mathfrak{m}-\sqrt{\mathfrak{m}^{2}-\Delta})$ is regarded as
the free energy of the relativistic Schr\"{o}dinger operator with
a particle of mass $\mathfrak{m}$; see the seminal paper of
Carmona {\it et al.} \cite{Carmona} for mathematical discussions
and its relation to L\'{e}vy processes. For  general $\alpha\in
(0,2)$, one may  refer to Ryznar \cite{Ryznar}, Baeumer {\it et
al.} \cite{Baeumer}, Kumara {\it et al.} \cite{Kumara}, and  the
references therein. RDEs have also played an essential role in the
theory of computer vision; see a special volume edited by Kimmel
{\it et al.} \cite{Kimmel}, in which  P.D.E. and scale-space
methods are focused and RDEs are particularly employed.
%
%

In this article, we consider the random initial data $u_{0}$ to be
subordinated Gaussian random fields and study the large-scale and
the small-scale  limits for the properly re-scaled solution field.
We prove that the two parameters $\alpha$  and $\mathfrak{m}> 0$
 play distinct roles in the two
scaling behaviors. For the large-scale limit (Theorem
\ref{thm:WeaklyDependenceMarco} and Theorem \ref{Macro Scaling
Theorem}), it is the mass $\mathfrak{m}> 0$ dominates the
space-time scaling and also the limiting field, which brings the
$\mathfrak{m}> 0$ in its structure. While for the small-scale
limit (Theorem \ref{thm:WeaklyDependenceMicro} and Theorem
\ref{Micro Scaling Theorem}), it is the spatial index $\alpha$
dominates both the scaling factor and the limiting field, and it
appears to be irrelevant for $\mathfrak{m}$ being positive or
zero.

In our discussions, the large-scale  Theorem
\ref{thm:WeaklyDependenceMarco} and Theorem \ref{Macro Scaling
Theorem} are respectively comparable  to the Central Limit Theorem
for local functionals of random fields with weak dependence
 in \cite{BreuerMajor}, and to a certain non-Gaussian Central Limit
Theorem for which the papers \cite{Taqqu,Dobrushin and Major} are
pioneering. For the small-scale Theorem
\ref{thm:WeaklyDependenceMicro} and Theorem \ref{Micro Scaling
Theorem}, they involve not only the space-time scaling on
$u(t,\mathbf{x}),\ t>0,\ \mathbf{x}\in \mathbb{R}^{n}$, but also
need to re-scale the initial data; to our knowledge, these are new
type results for the literature; see \cite{LSh2010} for the
authors' very recent study. As for the methodology for proofs, for
the Gaussian limits we employ the moments and the Feymann-type
diagrams used notably in \cite{BreuerMajor}, and for the
non-Gaussian limits we employ the truncation of Hermite expansions
used notably in \cite{AL,AnhHomo}.

We remark that, in  the non-relativistic case, i.e.
$\mathfrak{m}=0$, the large-scale limits for the random initial
value problem with multiple It\^{o}-Wiener integrals as input have
been discussed in Anh and Leonenko \cite{AL,AnhHomo}; subsequent
works, together with Burgers' equation, in this direction by the
authors and collaborators can be seen in
\cite{Higherorder,Barndorff,KLR,LM,LW,LW1,RAAfractionalBurger} and
the references therein. However, the multi-scaling limits due to
the different roles of the mass and the fractional-index, the
target of this article,  are {\em at all not} in the cited papers.
Moreover, in this article we are able to drop-off the usually
imposed isotropic assumption of the initial datum.

We should also mention that, in an article discussing tempered
stable L\'{e}vy processes by Rosi\'{n}ski \cite{Rons}, the author
proves rigorously, among others, the statement that such a process
in a short time looks a stable process while in a large time scale
it looks like a Brownian motion. This article has surfaced nicely
how the multi-scaling limits appear in the context of stochastic
processes (We are indebted  to the referee for indicating to us
the article \cite{Rons} and the relevant concept).


In Section 2, we  present some preliminaries; we state our main
results in Section 3, and all the proofs of our results are given
in Section 4.

Finally, we mention that the study on the PDEs with random initial
conditions can be traced back to  \cite{KF} and \cite{R}. Besides
the above mentioned literature, there also has very significant
progress on Burgers equation with different types of random input;
see the monograph of Woyczy\'{n}ski \cite{W lecture} and the
Chapter 6 of Bertoin \cite{B}.
%
%
%

{\bf Acknowledgement.} G.-R. Liu is  partially  supported by a
Taiwan NSC grant for graduate students. This article is mainly
worked while N.-R. Shieh visited York University (Canada) in Fall
2011 and Chinese University of Hong Kong in Spring 2012; the
hospitality and the financial support are appreciated. The content
of this article has been reported by N.-R. Shieh at the
probability scientific session of the Canadian Mathematical
Society 2011 Winter Meeting.

\vskip 20 pt

\section{Preliminaries}
\subsection{Green function for RDEs}
As understood,  we regard the spatial operator in the RDE
(\ref{relativistic diffusion equation}) as a psudo-differential
operator, see for example the book and the paper by Wong
\cite{Wong,Wong92}; the Green function, denoted by
$G_{\alpha,\mathfrak{m}}(t,\mathbf{x}),\ t>0, \mathbf{x}\in
\mathbb{R}^{n}$, for the Cauchy problem (\ref{relativistic
diffusion equation}) is thus determined by the (spatial) Fourier
transform $\widehat{G}_{\alpha,\mathfrak{m}}(t,\lambda),\
\alpha\in(0,2),\ \mathfrak{m}>0,$ which is given by
\begin{align}\label{fractional Green function}
\int_{\mathbb{R}^{n}}e^{i<\lambda,\mathbf{x}>}G_{\alpha,\mathfrak{m}}(t,\mathbf{x})d\mathbf{x}
=e^{-t\{(\mathfrak{m}^{\frac{2}{\alpha}}+|\lambda|^{2})^{\frac{\alpha}{2}}-\mathfrak{m}\}},\
\lambda\in \mathbb{R}^{n}.
\end{align}
See Carmona {\it et al.} \cite{Carmona} for $\alpha=1$ and  Ryznar
\cite{Ryznar} for  general $\alpha\in (0,2)$ (\cite{Ryznar} also
considers the boundary problem).  These papers also study
$G_{\alpha,\mathfrak{m}}(t,\mathbf{x}),\ \mathfrak{m}>0,$ as the
transition probability density of a L\`{e}vy process
$X_{\alpha,\mathfrak{m}}(t)$ which is the subordination of the
Brownian motion by a certain subordinator. The explicit expression
for the Green function is known only in the case  $\alpha=1$; see
for example the recent works of
 \cite{Baeumer,Kumara}, which give  explicit calculations  to show  that
the subordinator is normal inverse Gaussian.

The solution of (\ref{relativistic diffusion equation})
 is given in the  form
\begin{align}\label{meansquaresolution}
u(t,\mathbf{x};u_{0}(\cdot))=\int_{\mathbb{R}^{n}}
G_{\alpha,\mathfrak{m}}(t,\mathbf{x}-\mathbf{y})u_{0}(\mathbf{y})d\mathbf{y}.
\end{align}
In this work, our initial data is a second-order homogeneous
random field on $\mathbb{R}^n$, and thus the solution of
(\ref{relativistic diffusion equation}) should be understood as a
mean-square solution; resulting a spatial-temporal random solution
field $u(t,\mathbf{x})$; see  \cite[Proposition
1]{RAAfractionalBurger} for some discussions on the mean-square
solutions of parabolic PDEs with mean-square random initial data.


\subsection{Subordinated Gaussian fields as initial data
}

Let $(\Omega, \mathcal{F}, \mathcal{P})$  be an underlying
probability space,
 such that all random elements appeared in this article are measurable with respect to it.
We specify the initial data $u_{0}(\mathbf{x})$ be a subordinated
Gaussian field,  which is introduced by Dobrushin
\cite{Dobrushin1979}, as follows; see also \cite{AL,AnhHomo} for
more recent discussions.

\begin{bfseries}Condition A.\end{bfseries}
The initial data of  (\ref{relativistic diffusion equation}) is
assumed  to be  a random field on $\mathbb{R}^{n}$ given by
\begin{align}\label{initial data form}
u_{0}(\mathbf{x})=h(\zeta(\mathbf{x})),\ \mathbf{x}\in \mathbb{R}^{n},
\end{align}
where $\zeta(\mathbf{x})$ is a mean-square continuous and
homogeneous Gaussian random field with mean zero and variance 1,
and its spectral measure $F(d\lambda)$ has the (spectral) density
$f(\lambda) ,\ \lambda\in \mathbb{R}^{n}$; moreover,
$h:\mathbb{R}\rightarrow\mathbb{R}$ is a (non-random) function
such that
\begin{align}\label{L2 function}
\mathbb{E}h^{2}(\zeta(\mathbf{0}))=\int_{\mathbb{R}}h^{2}(r)p(r)dr<\infty;
\ \ \ p(r)=\frac{1}{\sqrt{2\pi}}e^{-\frac{r^{2}}{2}},\ r\in
\mathbb{R}.
\end{align}


Under Condition A,
by the Bochner-Khintchine theorem, we have the following spectral representation for the
covariance function of the Gaussian field $\zeta(\mathbf{x})$:
\begin{equation}\label{covariancespectralrep}
R(\mathbf{x})=\textup{Cov}(\zeta(\mathbf{0}),\zeta(\mathbf{x}))=\int_{\mathbb{R}^{n}}e^{i<\lambda,\mathbf{x}>}f(\lambda)d\lambda.
\end{equation}
Moreover, by the Karhunen Theorem, $\zeta(\mathbf{x})$ has the
representation
\begin{align}\label{sample path represent}
\zeta(\mathbf{x})=\int_{\mathbb{R}^{n}}
e^{i<\lambda,\mathbf{x}>}\sqrt{f(\lambda)}W(d\lambda),\ \mathbf{x}\in \mathbb{R}^{n},
\end{align}
where $W(d\lambda)$ is the standard complex-valued Gaussian white
noise on the Fourier domain $\mathbb{R}^{n}$; that is, a centered
orthogonal-scattered Gaussian random  measure on $\mathbb{R}^{n}$
such that
 $W(\Delta_{1})=\overline{W(-\Delta_{1})}$
and $\mathbb{E}W(\Delta_{1})
\overline{W(\Delta_{2})}=\textup{Leb}(\Delta_{1}\cap\Delta_{2})$
for any
 $\Delta_{1},\Delta_{2}\in
\mathcal{B}(\mathbb{R}^{n})$. See, for example, the book of
Leonenko \cite[Theorem 1.1.3]{L} for the above facts.
We need the following expansion of $h(r)$ in the Hilbert space
$L^{2}(\mathbb{R},p(r)dr)$:
\begin{align}\label{hermiteexpansion}
h(r)=C_{0}+\sum_{l=1}^{\infty}{C_{l}
\frac{H_{l}(r)}{\sqrt{l !}}},
\end{align}
where
\begin{align}\label{hermitecoeff}
C_{l}=\int_{\mathbb{R}}h(r)\frac{H_{l}(r)}{\sqrt{l!}}p(r)dr,
\end{align}
and $\{H_{l}(r),\ l=0,1,2,...\}$ are the Hermite
polynomials, that is,
\[H_{l}(r)=(-1)^{l}e^{\frac{r^{2}}{2}}\frac{d^{l}}{dr^{l}}e^{-\frac{r^{2}}{2}},
\ \ \textup{for}\ \ l\in \{0,1,2,...\}.\]
Accordingly, the {\it Hermite rank} of the function $h(\cdot)$ is
defined by
\begin{align*}
m:=\textup{inf}\{l\geq 1:\ C_{l}\neq 0\}.
\end{align*}
%
It is well-known that(see, for example, Major \cite[Corollary 5.5
and p. 30]{M}):
\begin{align}\label{expectionhermite}
\mathbb{E}[H_{l_{1}}(\zeta(\mathbf{y}))H_{l_{2}}(\zeta(\mathbf{z}))]=
\delta^{l_{1}}_{l_{2}} l_{1}!
R^{l_{1}}(\mathbf{y}-\mathbf{z}),\ \ \ \mathbf{y},\ \mathbf{z}\in
\mathbb{R}^{n},
\end{align}
($\delta^{\sigma_{1}}_{\sigma_{2}}$ is the Kronecker symbol) and
\begin{align}\label{itoformula}
H_{l}(\zeta(\mathbf{x}))=
\int^{'}_{\mathbb{R}^{n\times l}}e^{i<\mathbf{x},\lambda_{1}+...+\lambda_{l}>}\prod_{k=1}^{l}{\sqrt{f(\lambda_{k})}}
W(d\lambda_{k}).
\end{align}
In the above, (\ref{itoformula}) means the   {\it multiple
It\^{o}-Wiener integral} representation  and the integration
$\int^{'}$ means that it excludes the diagonal hyperplanes
$\mathbf{z}_{i}=\mp \mathbf{z}_{j},\ i, j=1,...,l, i\neq j$.

We impose two different conditions on the singularity  of the
spectral density $f(\lambda)$ at $\mathbf{0}$, which yield,
respectively, the Gaussian and the non-Gaussian scaling-limits.


\begin{bfseries}Condition B.\end{bfseries}
The spectral density function
$f(\lambda)$
of the Gaussian random field $\zeta(\mathbf{x})$ in Condition A
can be written as
\begin{equation}
f(\lambda)=\frac{B(\lambda)}{|\lambda|^{n-\kappa}}\ \textup{for some}\ \kappa>\frac{n}{m},
\end{equation}
where $m$ is the Hermite rank of the  function $h$, and the
$B(\cdot)\in\mathrm{C}(\mathbb{R}^{n})$ is of suitable decay at
infinity to ensure   $f\in L^{1}(\mathbb{R}^{n})$.

\begin{bfseries}Condition C.\end{bfseries}
The spectral density function $f(\lambda)$ of the Gaussian random
field $\zeta(\mathbf{x})$ in Condition A can be written as
\begin{align}\label{long}
f(\lambda)= \frac{B(\lambda)}{|\lambda|^{n-\kappa}},\ \
0<\kappa<\frac{n}{m},
\end{align}
where $m$ is the Hermite rank of the function $h$, and the
$B(\cdot)\in \mathrm{C}(\mathbb{R}^{n})$ is of  suitable decay at
infinity to ensure $f\in L^{1}(\mathbb{R}^{n})$, and moreover
$B(\mathbf{0})>0$.
%

Note that,  in the two conditions, we do not assume that the
$B(\cdot)$ is radial in $\cdot$, so that the field
$u_0(\mathbf{x})$ is not necessary to be isotropic. We also
mention that, the Condition B means that the density $f$ either is
regular at $\mathbf{0}$, or has  a singularity for which the order
is less than $n(1-1/m)$; while the Condition C means that $f$ has
a singularity at $\mathbf{0}$ for which the order is higher than
$n(1-1/m)$.

By (\ref{covariancespectralrep}) and the convolutions, we have,
for each $l\ge 1$,
\begin{equation}\label{powercovariancespectralrep}
R^{l}(\mathbf{x})=\int_{\mathbb{R}^{n}}e^{i<\lambda,\mathbf{x}>}f^{*l}(\lambda)d\lambda,\ l\in \mathbb{N},
\end{equation}
where $f^{*l}(\lambda)$ is the $l$-fold  convolution of $f$
defined recursively as: $f^{*1}=f$ and
\begin{equation*}
f^{*l}(\lambda)=\int_{\mathbb{R}^{n}}f(\lambda-\eta)f^{*(l-1)}(\eta)d\eta,\
 l\ge 2.
\end{equation*}
The following analytic lemma asserts  the behavior of $f^{*l},\
l\in \mathbb{N};$ for completeness, we give its  proof  in
Appendix A.
\begin{Lemma}
Suppose that the spectral density function $f$ has the form,

\begin{equation*}
f(\lambda)=\frac{B(\lambda)}{|\lambda|^{n-\kappa}},\ \kappa>0,
\end{equation*}
for some non-negative bounded and continuous function $B(\lambda)$
so that  $f\in L^{1}(\mathbb{R}^{n})$. Then for any $k\geq 2$
there exists a bounded function $B_{k}\in
\mathrm{C}(\mathbb{R}^{n}\backslash\{\mathbf{0}\})$ such that the
$k$-fold convolution $f^{*k}$ of $f$ can be written as
\begin{equation}\label{singularconvolution1}
f^{*k}(\lambda)=
\left\{\begin{array}{lr}
B_{k}(\lambda)|\lambda|^{k\kappa-n},\ \ & \textup{for}\ k\kappa<n,
\\
B_{k}(\lambda)\textup{ln}(2+\frac{1}{|\lambda|}),\ \ & \textup{for}\ k\kappa=n,
\\
B_{k}(\lambda)\in \mathrm{C}(\mathbb{R}^{n}),\ \ & \textup{for}\ k\kappa> n.
\end{array}
\right.
\end{equation}
Moreover, for any $k_{1}>k_{2}>n/\kappa$ the inequality
$\underset{\lambda\in \mathbb{R}^{n}}{\textup{sup}}B_{k_{1}}(\lambda)\leq
\underset{\lambda\in \mathbb{R}^{n}}{\textup{sup}}B_{k_{2}}(\lambda)
$ holds.
\end{Lemma}

To understand the difference of the Conditions B and C, in view of
Lemma 1, the Condition B implies that the $k$-fold convolution
$f^{*k},$ $k\geq m$, has no singularity at the origin
$\lambda=\mathbf{0}$, which in turn assets that  the spectral
density of the random initial data $u_{0}$ has no singularity at
$\lambda=\mathbf{0}$; while the Condition C asserts that the
initial data $u_0$ has a spectral density which is singularity at
$\lambda=\mathbf{0}$. The situation  can be described as,
respectively, the long-range and  the short-range dependence of
the initial field $u_0$; a central notion in vast applications, as
one may refer to the special volume by Doukhan, Oppenheim, and
Taqqu \cite{DOT}.


\section{Main results}

The significant difference between the Condition B and the
Condition C, as remarked at the end of the last section, is
employed to obtain the Gaussian and respectively the non-Gaussian
scaling-limits. We will present them in two subsections.

In the context henceforth, the notation $\Rightarrow$ denotes the
convergence of random variables (respectively,  random families)
in the sense of distribution (respectively, finite-dimensional
distributions).

\subsection{Gaussian limits with initial data in (A,B)}

As mentioned in the Section 1, we will present the large-scale and
the small-scale limit theorems. We remark that  our Theorems 1 and
2 in this subsection are comparable to  the central limit theorem
for local functionals of random fields with weak dependence in
Breuer and Major \cite{BreuerMajor}. The novel feature is that the
the mass $\mathfrak{m}>0$ and the fractional-index $\alpha$ play
different roles in the two-scales.

\begin{Theorem}\label{thm:WeaklyDependenceMarco}
Let $u(t,\mathbf{x};u_{0}(\cdot)),\ t>0,\ \mathbf{x}\in \mathbb{R}^{n}$, be the mean-square solution
of (\ref{relativistic diffusion equation}) with $\mathfrak{m}>0$
and the initial data $u_{0}(\mathbf{x})=h(\zeta(\mathbf{x}))$
satisfy Condition A and B with the Hermite rank $m\geq 1$.
Then when $T\rightarrow\infty$,
\begin{equation*}
T^{\frac{n}{4}}\Big\{ u(T t,\sqrt{T}\mathbf{x};u_{0}(\cdot))-C_{0}
\Big\} \Rightarrow U(t,\mathbf{x}),
\end{equation*}
where
$U(t,\mathbf{x}),\ t>0,\ \mathbf{x}\in \mathbb{R}^{n}$,
is a Gaussian field
with the following spectral representation:
\begin{equation}\label{thm:WeaklyDependenceMarcollimitingfield}
U(t,\mathbf{x})=
\int_{\mathbb{R}^{n}}
e^{i<\lambda,\mathbf{x}>}
\sigma_{m}e^{-t\frac{\alpha}{2}\mathfrak{m}^{1-\frac{2}{\alpha}}|\lambda|^{2}}W(d\lambda),\
\sigma_{m}=\Big(\overset{\infty}{\underset{r=m}{\sum}}
f^{*r}(\mathbf{0})
C^{2}_{r}\Big)^{\frac{1}{2}},
\end{equation}
where $W(d\lambda)$
is a complex-valued standard Gaussian noise measure on $\mathbb{R}^{n}$ (c.f. (\ref{sample path represent})).
\end{Theorem}

For the small-scale limit, we need to re-scale the initial data
too; thus the notation
$u_{0}(\varepsilon^{-\frac{1}{\alpha}-\chi}\cdot)$ imposed on
$u_{0}$ wants to mean that the variable of $u_{0}$ is under the
indicated dilation factor $\varepsilon^{-\frac{1}{\alpha}-\chi}$.

\begin{Theorem}\label{thm:WeaklyDependenceMicro}
Let $u(t,\mathbf{x};u_{0}(\cdot)),\ t>0,\ \mathbf{x}\in \mathbb{R}^{n}$, be the mean-square solution
of (\ref{relativistic diffusion equation}) with $\mathfrak{m}>0$
and the initial data $u_{0}(\mathbf{x})=h(\zeta(\mathbf{x}))$
satisfy Condition A and B with the Hermite rank $m\geq 1$.
For any $\chi>0$, when $\varepsilon\rightarrow 0$,
\begin{equation}
\varepsilon^{-\frac{n\chi}{2}}\Big\{ u(\varepsilon
t,\varepsilon^{\frac{1}{\alpha}}\mathbf{x};u_{0}(\varepsilon^{-\frac{1}{\alpha}-\chi}\cdot))-C_0\Big\}
\Rightarrow V(t,\mathbf{x}),
\end{equation}
where $V(t,\mathbf{x}),\ t>0,\ \mathbf{x}\in \mathbb{R}^{n}$,
is a Gaussian field
with the following spectral representation:
\begin{equation}\label{thm:WeaklyDependenceMicrollimitingfield}
V(t,\mathbf{x})=
\int_{\mathbb{R}^{n}}
e^{i<\lambda,\mathbf{x}>}
\sigma_{m}e^{-t|\lambda|^{\alpha}}W(d\lambda),\
\sigma_{m}=\Big(\overset{\infty}{\underset{r=m}{\sum}}
f^{*r}(\mathbf{0})
C^{2}_{r}\Big)^{\frac{1}{2}},
\end{equation}
where $W(d\lambda)$
is a complex-valued standard Gaussian noise measure on $\mathbb{R}^{n}$.
\end{Theorem}

{\bf Remark.} The typical case for Theorem 2 is $\alpha=1,
\chi=1/2.$ In  this critical case,  the scaling order for Theorems
1 and 2 is the same, namely $n/4$. However,  the spatial scaling
is square-root in Theorem 1 while is linear in Theorem 2;
moreover, the integral kernel for  the limiting field in two
theorems is Gauss vs. Poisson. The latter situation can be
conferred to an analytic discussion in Wong \cite{Wong92}.

\subsection{Non-Gaussian limits  with  initial data  in (A,C)}

As in the above subsection, we have the large-scale and the
small-scale limits;  however the high singularity order in the
Condition C assets that  our limiting fields are now non-Gaussian.
The non-Gaussian limits of the convolution type; which can be seen
in the pioneering papers of Taqqu \cite{Taqqu} and Dobrushin and
Major \cite{Dobrushin and Major}, and  more recent Anh and
Leonenko \cite{AL,AnhHomo}.

\begin{Theorem}\label{Macro Scaling Theorem}
Let $u(t,\mathbf{x};u_{0}(\cdot)),\ t>0,\ \mathbf{x}\in
\mathbb{R}^{n},$ be the mean-square solution of (\ref{relativistic
diffusion equation}) whose initial data
$\{u_{0}(\mathbf{x})=h(\zeta(\mathbf{x})),\ \mathbf{x}\in
\mathbb{R}^{n}\}$ satisfy Condition A and C with $\kappa\in(0,
\frac{n}{m})$ and $1<m$, where $m$ is the Hermite rank of the
non-random function $h$ on $\mathbb{R}$, which has the Hermite
coefficients $C_{j},j=0,1,\ldots$. Then when $T\rightarrow
\infty$,
\begin{align}\label{rescaled macro field}
T^{\frac{m\kappa}{4}}
\Big\{u(Tt, \sqrt{T}\mathbf{x}
;h(\zeta(\cdot)))-C_{0}
\Big\}\Rightarrow U_{m}(t,\mathbf{x}),
\end{align}
where
$U_{m}(t,\mathbf{x})$ is represented by the following multiple Wiener integrals
\begin{align}\label{limiting field(macro)}
U_{m}(t,\mathbf{x})\hspace{-0.1cm}=\hspace{-0.1cm}B^{\frac{m}{2}}(\mathbf{0})\frac{C_{m}}{\sqrt{m!}}
\int_{\mathbb{R}^{n\times m}}^{'}\hspace{-0.8cm}
e^{i<\mathbf{x},\lambda_{1}+\cdots+\lambda_{m}>}
\frac{
\textup{exp}(-t\frac{\alpha}{2}\mathfrak{m}^{1-\frac{2}{\alpha}}|\lambda_{1}+\cdots+\lambda_{m}|^{2})
} {(|\lambda_{1}| \cdots
|\lambda_{m}|)^{\frac{n-\kappa}{2}}}\overset{m}{\underset{l=1}{\prod}}W(d\lambda_{l}),
\end{align}
where $\int^{'}_{\mathbb{R}^{n\times m}}\cdots$ denotes an
$m$-fold Wiener integral with respect to the complex Gaussian
white noise $W(\cdot)$ on $\mathbb{R}^{n}$.
\end{Theorem}

\begin{Theorem}\label{Micro Scaling Theorem}
Let $u(t,\mathbf{x};u_{0}(\cdot))$ be the mean-square solution to
(\ref{relativistic diffusion equation}) whose initial data
$\{u_{0}(\mathbf{x})=h(\zeta(\mathbf{x})),\ \mathbf{x}\in
\mathbb{R}^{n}\}$ satisfy Condition A and C with $\kappa\in(0,
\frac{n}{m})$ and $1<m$, where $m$ is the Hermite rank of the
function $h$. Then, for any fixed parameter $\chi>0$, when
$\varepsilon\rightarrow 0$,
\begin{equation}\label{rescaled micro field}
\varepsilon^{-\frac{m\kappa\chi}{2}}
\Big\{u(\varepsilon t, \varepsilon^{\frac{1}{\alpha}}\mathbf{x}
;h(\zeta((\varepsilon^{-\frac{1}{\alpha}-\chi})\cdot)))-C_{0}
\Big\}\Rightarrow V_{m}(t,\mathbf{x}),
\end{equation}
where
$V_{m}(t,\mathbf{x})$ is represented by the multiple Wiener integrals
\begin{align}\label{limiting field(micro)}
V_{m}(t,\mathbf{x})\hspace{-0.1cm}=\hspace{-0.1cm}
B^{\frac{m}{2}}(\mathbf{0})
\frac{C_{m}}{\sqrt{m!}}
\int_{\mathbb{R}^{n\times m}}^{'}\hspace{-0.6cm}
e^{i<\mathbf{x},\lambda_{1}+\cdots+\lambda_{m}>}
\frac{
\textup{exp}(-
t|\lambda_{1}+\ldots+\lambda_{m}|^{\alpha})
} {(|\lambda_{1}| \cdots
|\lambda_{m}|)^{\frac{n-\kappa}{2}}}\overset{m}{\underset{l=1}{\prod}}W(d\lambda_{l}).
\end{align}
\end{Theorem}

{\bf Remark.} In  \cite{AnhHomo} the authors considered a hybrid
differential operator in the spatial variable (the Riesz-Bessel
operator), as follows
$$
-(-\Delta)^{\alpha/2}(I-\Delta)^{\gamma/2},\ \alpha\in(0,2),\
\gamma\ge0.
$$
{\em However}, in their main Theorem 2.3, a large-scale limit in
our context,  only the Riesz parameter $\alpha$ plays the role and
the Bessel parameter $\gamma$ is invisible. This intrigue
situation is now justified by the RFD (\ref{relativistic diffusion
equation}), which we could say that it is  ``physically correct"
to consider the relativistic operator
$(\mathfrak{m}-(\mathfrak{m}^{\frac{2}{\alpha}}-\Delta)^{\frac{\alpha}{2}})$
rather than the Bessel operator in the form presented in
\cite{AnhHomo}.

\vskip 20 pt

\section{Proofs of Theorems}

The following two-scale property of the relativistic Green
function $G_{\alpha,\mathfrak{m}}$ is the key to our results; when
one deals the Laplacian or the fractional-Laplacian operator, it
is instead only the mono-scaling. We describe this two-scale
property in terms of Fourier transforms.

\begin{equation}\label{proofmacroexponent}
\widehat{G}_{\alpha,\mathfrak{m}}(Tt,T^{-\frac{1}{2}}\lambda)=
\textup{exp}\Big\{Tt(\mathfrak{m}-(\mathfrak{m}^{\frac{2}{\alpha}}+T^{-1}|\lambda|^{2})^{\frac{\alpha}{2}})\Big\}
\rightarrow
\textup{exp}\Big\{-t\frac{\alpha}{2}\mathfrak{m}^{1-\frac{2}{\alpha}}|\lambda|^{2}\Big\},
\end{equation}
as  $T\rightarrow\infty$; (\ref{proofmacroexponent}) is a
consequence of the of Taylor's expansion,
\begin{align*}
\mathfrak{m}-(\mathfrak{m}^{\frac{2}{\alpha}}+T^{-1}|\lambda|^{2})^{\frac{\alpha}{2}}
=&
\mathfrak{m}-
\Big(\mathfrak{m}+\frac{\alpha}{2}(\mathfrak{m}^{\frac{2}{\alpha}})^{\frac{\alpha}{2}-1}T^{-1}|\lambda|^{2}
+\frac{\alpha}{4}(\frac{\alpha}{2}-1)c_{T}^{\frac{\alpha}{2}-2}T^{-2}|\lambda|^{4}\Big)
\\=&
-\frac{\alpha}{2}(\mathfrak{m}^{\frac{2}{\alpha}})^{\frac{\alpha}{2}-1}T^{-1}|\lambda|^{2}
+\frac{\alpha}{4}(1-\frac{\alpha}{2})c_{T}^{\frac{\alpha}{2}-2}T^{-2}|\lambda|^{4}
\end{align*}
for some $c_{T}\in (\mathfrak{m}^{\frac{2}{\alpha}},\mathfrak{m}^{\frac{2}{\alpha}}+T^{-1}|\lambda|^{2})$.
In contrast to the large-scale (\ref{proofmacroexponent}), we have
the following small-scale, as  $\varepsilon\rightarrow0$,
\begin{equation}\label{rescaledgreenlimit}
\widehat{G}_{\alpha,\mathfrak{m}}(\varepsilon
t,\varepsilon^{-\frac{1}{\alpha}}\lambda) = e^{\varepsilon
t\mathfrak{m}}e^{-\varepsilon t
(\mathfrak{m}^{\frac{2}{\alpha}}+\varepsilon^{-\frac{2}{\alpha}}|\lambda|^{2})^{\frac{\alpha}{2}}}
\rightarrow e^{-t|\lambda|^{\alpha}}.
\end{equation}
We observe that (\ref{rescaledgreenlimit}) indeed holds no matter $\mathfrak{m}$ is $>0$ or $=0$.%
%

\bigskip

\noindent {\bf Proofs of Theorems 1 and 2.}

\noindent We apply the Hermite expansion (\ref{hermiteexpansion})
to $u(t,x)$, For the large-scale, we set
\begin{align*}
&X_{T}(t,\mathbf{x})=T^{\frac{n}{4}}u(T t,\sqrt{T}\mathbf{x};u_{0}(\cdot))-C_{0}
\\
=& T^{\frac{n}{4}}
\overset{\infty}{\underset{k=m}{\sum}}\frac{C_{k}}{\sqrt{k!}}\int_{\mathbb{R}^{n}}
G_{\alpha,\mathfrak{m}}(T t,\sqrt{T}\mathbf{x}-\mathbf{y})
H_{k}(\zeta(\mathbf{y}))d\mathbf{y},
\end{align*}
and for the small-scale, we set
\begin{align*}
Y_{\varepsilon}(t,\mathbf{x}):=&\varepsilon^{-\frac{n\chi}{2}}u(\varepsilon
t,\varepsilon^{\frac{1}{\alpha}}\mathbf{x};u_{0}(\varepsilon^{-\frac{1}{\alpha}-\chi}\cdot))-C_{0}
\\
=& \varepsilon^{-\frac{n\chi}{2}}
\overset{\infty}{\underset{l=m}{\sum}}\frac{C_{l}}{\sqrt{l!}}\int_{\mathbb{R}^{n}}
G_{\alpha,\mathfrak{m}}(\varepsilon
t,\varepsilon^{\frac{1}{\alpha}}\mathbf{x}-\mathbf{y})
H_{l}(\zeta(\varepsilon^{-\frac{1}{\alpha}-\chi}\mathbf{y}))d\mathbf{y}.
\end{align*}

In the below, we only proceed the proof of Theorem 2, the
small-scale limit, and see how the rescaling of the initial data
is needed to obtain the desired limit; the proof of Theorem 1 is
parallel, and does not require the rescaling of the initial data.
Since the proof in the following   does not require the
$\mathfrak{m}$ to be strictly positive, our Theorem 2 also
provides a small-scale version of the large-scale, i.e. the usual,
limit result in \cite{AnhHomo}. The methodology of the proof can
be traced back to \cite{BreuerMajor}.

For any $M\in \mathbb{N}$ and any set of real numbers
$\{a_{1},a_{2},\cdots,a_{M}\}$, denote
\begin{align}
\xi_{\varepsilon}:=
\overset{M}{\underset{j=1}{\sum}}a_{j}Y_{\varepsilon}(t_{j},\mathbf{x}_{j}),
\end{align}
where $\{t_{1},\cdots,t_{M}\}\subset \mathbb{R}_{+}$ and
$\{\mathbf{x}_{1},\cdots,\mathbf{x}_{M}\}\subset \mathbb{R}^{n}$
are arbitrary. In order to apply the Method of Moments to prove
the statement of Theorem \ref{thm:WeaklyDependenceMicro}, we need
to verify the following:
\begin{align}\label{Markovmethod}
\underset{\varepsilon\rightarrow 0}{\textup{lim}}\ \mathbb{E}\xi_{\varepsilon}^{p}=
\left
\{\begin{array}{lr}0, & p=2\nu+1
\\ (p-1)!!
\Big\{\mathbb{E}\Big[\big(\overset{M}{\underset{j=1}{\sum}}a_{j}V(t_{j},\mathbf{x}_{j})\big)^{2}\Big]\Big\}^{\nu},
& p=2\nu
\end{array}
\right.,
\end{align}
where $V(t,\mathbf{x})$ is defined in (\ref{thm:WeaklyDependenceMicrollimitingfield}).
%
%
We remark that the high (i.e. $p>2$) moments is needed, since
$\xi_{\varepsilon}$ is not Gaussian, though the wanted limit is
Gaussian. Firstly, we split $\xi_{\varepsilon}$ into two parts:
\begin{equation}\label{proofthm:weaksmall3}
\xi_{\varepsilon}=\xi_{\varepsilon,\leq N}+\xi_{\varepsilon,>N},
\end{equation}
where (henceforth, we will suppress the indices $\alpha$ and
$\mathfrak{m}$ for $G_{\alpha,\mathfrak{m}}$ and
$\widehat{G}_{\alpha,\mathfrak{m}}$)
\begin{align}
\xi_{\varepsilon,>N}=\overset{M}{\underset{j=1}{\sum}}a_{j}
\varepsilon^{-\frac{n\chi}{2}}
\overset{\infty}{\underset{l=N+1}{\sum}}\frac{C_{l}}{\sqrt{l!}}\int_{\mathbb{R}^{n}}
G(\varepsilon
t_{j},\varepsilon^{\frac{1}{\alpha}}\mathbf{x}_{j}-\mathbf{y})
H_{l}(\zeta(\varepsilon^{-\frac{1}{\alpha}-\chi}\mathbf{y}))d\mathbf{y},
\end{align}
and we prove that $E[\xi_{\varepsilon,>N}^2]\rightarrow 0$,
whenever $N$ is chosen large enough. Observe that for any $N\geq
m-1$, by (\ref{expectionhermite}),
\begin{align}\notag
&\mathbb{E}(\xi_{\varepsilon,>N})^{2}=\mathbb{E}\Big[\big(\overset{M}{\underset{j=1}{\sum}}a_{j}
\varepsilon^{-\frac{n\chi}{2}}
\overset{\infty}{\underset{l=N+1}{\sum}}\frac{C_{l}}{\sqrt{l!}}\int_{\mathbb{R}^{n}}
G(\varepsilon t_{j},\varepsilon^{\frac{1}{\alpha}}\mathbf{x}_{j}-\mathbf{y})
H_{l}(\zeta(\varepsilon^{-\frac{1}{\alpha}-\chi}\mathbf{y}))d\mathbf{y}
\big)^{2}\Big]
\\\notag
=&
\overset{M}{\underset{j_{1},j_{2}=1}{\sum}}a_{j_{1}}a_{j_{2}}
\varepsilon^{-n\chi}
\overset{\infty}{\underset{l=N+1}{\sum}}C^{2}_{l}\int_{\mathbb{R}^{2n}}
\hspace{-0.16cm}G(\varepsilon t_{j_{1}},\varepsilon^{\frac{1}{\alpha}}\mathbf{x}_{j_{1}}-\mathbf{y}_{1})
G(\varepsilon t_{j_{2}},\varepsilon^{\frac{1}{\alpha}}\mathbf{x}_{j_{2}}-\mathbf{y}_{2})
R^{l}(\varepsilon^{-\frac{1}{\alpha}-\chi}(\mathbf{y}_{1}-\mathbf{y}_{2})))
\\\label{proofthm:weaksmall1}
=&
\overset{M}{\underset{j_{1},j_{2}=1}{\sum}}a_{j_{1}}a_{j_{2}}
\varepsilon^{-n\chi}
\overset{\infty}{\underset{l=N+1}{\sum}}C^{2}_{l}\int_{\mathbb{R}^{n}}
G(\varepsilon (t_{j_{1}}+t_{j_{2}}),\varepsilon^{\frac{1}{\alpha}}(\mathbf{x}_{j_{1}}-\mathbf{x}_{j_{2}})-\mathbf{z})
R^{l}(\varepsilon^{-\frac{1}{\alpha}-\chi}\mathbf{z})
d\mathbf{z},
\end{align}
where the last equality is followed by changing of variables, the symmetry property $G(t,\mathbf{z})=G(t,-\mathbf{z})$ of
the transition probability density function $G$, and its semigroup property
\begin{align*}
\int_{\mathbb{R}^{n}}
\hspace{-0.1cm}G(\varepsilon t_{j_{1}},\varepsilon^{\frac{1}{\alpha}}\mathbf{x}_{j_{1}}-(\mathbf{z}-\mathbf{z}^{'}))
G(\varepsilon t_{j_{2}},\varepsilon^{\frac{1}{\alpha}}\mathbf{x}_{j_{2}}-\mathbf{z}^{'})
d\mathbf{z}{'}
=
G(\varepsilon (t_{j_{1}}+t_{j_{2}}),\varepsilon^{\frac{1}{\alpha}}(\mathbf{x}_{j_{1}}-\mathbf{x}_{j_{2}})-\mathbf{z}).
\end{align*}
Continue to (\ref{proofthm:weaksmall1}),
by the spectral representation (\ref{powercovariancespectralrep}) for the $k$-th power of the covariance function $R(\cdot)$,
it is equal to
\begin{align}\notag
&
\overset{M}{\underset{j_{1},j_{2}=1}{\sum}}a_{j_{1}}a_{j_{2}}
\varepsilon^{-n\chi}
\overset{\infty}{\underset{l=N+1}{\sum}}C^{2}_{l}
\underset{\mathbb{R}^{n}}{\int}
G(\varepsilon (t_{j_{1}}+t_{j_{2}}),\varepsilon^{\frac{1}{\alpha}}(\mathbf{x}_{j_{1}}-\mathbf{x}_{j_{2}})-\mathbf{z})
\underset{\mathbb{R}^{n}}{\int}e^{i<\varepsilon^{-\frac{1}{\alpha}-\chi}\mathbf{z},\lambda>}f^{*l}(\lambda)d\lambda
d\mathbf{z}
\\\notag
=&
\overset{M}{\underset{j_{1},j_{2}=1}{\sum}}a_{j_{1}}a_{j_{2}}
\varepsilon^{-n\chi}
\overset{\infty}{\underset{l=N+1}{\sum}}C^{2}_{l}
\underset{\mathbb{R}^{n}}{\int}
e^{i\varepsilon^{-\chi}<\lambda,\mathbf{x}_{j_{1}}-\mathbf{x}_{j_{2}}>}
\widehat{G}(\varepsilon(t_{j_{1}}+t_{j_{2}}),\varepsilon^{-\frac{1}{\alpha}-\chi}\lambda)
f^{*l}(\lambda)d\lambda
\\\notag
=&
\overset{M}{\underset{j_{1},j_{2}=1}{\sum}}a_{j_{1}}a_{j_{2}}
\overset{\infty}{\underset{l=N+1}{\sum}}C^{2}_{l}
\underset{\mathbb{R}^{n}}{\int}
e^{i<\lambda,\mathbf{x}_{j_{1}}-\mathbf{x}_{j_{2}}>}
\textup{exp}\{\varepsilon(t_{j_{1}}+t_{j_{2}})
[\mathfrak{m}-(\mathfrak{m}^{\frac{2}{\alpha}}+|\varepsilon^{-\frac{1}{\alpha}}\lambda|^{2})^{\frac{\alpha}{2}}]\}
f^{*l}(\varepsilon^{\chi}\lambda)d\lambda
\\\label{proofthm:weaksmall2}
\rightarrow &
\overset{M}{\underset{j_{1},j_{2}=1}{\sum}}a_{j_{1}}a_{j_{2}}
\overset{\infty}{\underset{l=N+1}{\sum}}C^{2}_{l}f^{*l}(\mathbf{0})\int_{\mathbb{R}^{n}}
e^{i<\lambda,\mathbf{x}_{j_{1}}-\mathbf{x}_{j_{2}}>}
\textup{exp}\{-(t_{j_{1}}+t_{j_{2}})
|\lambda|^{\alpha}\}d\lambda<\infty
\end{align}
when $\varepsilon\rightarrow 0$, where $f^{*l}(\cdot),\ l\geq m$,
are continuous and uniformly bounded on $\mathbb{R}^{n}$
since
Condition B and Lemma 1 imply:
\begin{equation}\label{inequalityconvolutionweakdepsmall}
f^{*l}(\lambda)\hspace{-0.1cm}=\hspace{-0.1cm}\int_{\mathbb{R}^{n}}\hspace{-0.2cm}
f^{*m}(\lambda-\eta)f^{*(l-m)}(\eta)d\eta
\leq \parallel\hspace{-0.1cm} B_{m}\hspace{-0.1cm}\parallel_{\infty}
\hspace{-0.15cm}\int_{\mathbb{R}^{n}}\hspace{-0.2cm}f^{*(l-m)}(\eta)d\eta
=\parallel\hspace{-0.1cm}B_{m}\hspace{-0.1cm}\parallel_{\infty}\ \ \forall l>m.
\end{equation}
%
From (\ref{proofthm:weaksmall2}), for any $\delta>0$ there exists $N_{0}\in \mathbb{N}$,
$\varepsilon_{0}>0$
such that
\begin{align}\label{proofthm:weaksmall4}
\mathbb{E}(\xi_{\varepsilon,>N})^{2}<\delta,\ \textup{for any}\ N\geq N_{0},\ \varepsilon<\varepsilon_{0},
\end{align}
which implies that we suffice to prove a truncated  version of
(\ref{Markovmethod}) as follows:
\begin{align}\label{Markovmethod2}
\underset{\varepsilon\rightarrow 0}{\textup{lim}}\ \mathbb{E}\xi_{\varepsilon,\leq N_{0}}^{p}=
\left
\{\begin{array}{lr}0, & p=2\nu+1
\\ (p-1)!!
\Big\{\mathbb{E}\Big[\big(\overset{M}{\underset{j=1}{\sum}}a_{j}V_{m,N_{0}}(t_{j},\mathbf{x}_{j})\big)^{2}\Big]\Big\}^{\nu},
& p=2\nu
\end{array}
\right.,
\end{align}
where
\begin{align}
V_{m,N_{0}}(t,\mathbf{x})=
\int_{\mathbb{R}^{n}}
e^{i<\lambda,\mathbf{x}>}
\sigma_{m,N_{0}}e^{-t|\lambda|^{\alpha}}W(d\lambda)\ \ \textup{with}\
\sigma_{m,N_{0}}=(\overset{N_{0}}{\underset{r=m}{\sum}}
f^{*r}(\mathbf{0})
C^{2}_{r})^{\frac{1}{2}}.
\end{align}
By (\ref{proofthm:weaksmall3}) for the definition of
$\xi_{\varepsilon,\leq
N_{0}}(=\xi_{\varepsilon}-\xi_{\varepsilon,>N_{0}})$, and our
rescaling of the initial data, we have
\begin{align}\notag
\mathbb{E}(\xi_{\varepsilon,\leq N_{0}})^{p}=&\varepsilon^{-\frac{pn\chi}{2}}
\overset{M}{\underset{{j_{1},\cdots,j_{p}=1}}{\sum}}\
\overset{N_{0}}{\underset{l_{1},\cdots,l_{p}=m}{\sum}}
\Big[\overset{p}{\underset{i=1}{\prod}}a_{j_{i}} \frac{C_{l_{i}}}{\sqrt{l_{i}!}}
\Big]
\\\notag
&\times
\int_{\mathbb{R}^{np}}
\Big\{
\overset{p}{\underset{i=1}{\prod}}
G(\varepsilon t_{j_{i}},\varepsilon^{\frac{1}{\alpha}}\mathbf{x}_{j_{i}}-\mathbf{y}_{i})\Big\}
\Big[\mathbb{E}\overset{p}{\underset{i=1}{\prod}} H_{l_{i}}(\zeta(\varepsilon^{-\frac{1}{\alpha}-\chi}\mathbf{y}_{i}))\Big]
d\mathbf{y}_{1}\cdots d\mathbf{y}_{p}
\\\label{proofthm:weaksmall5}
=
&\varepsilon^{-\frac{pn\chi}{2}}
\overset{M}{\underset{{j_{1},\cdots,j_{p}=1}}{\sum}}\
\overset{N_{0}}{\underset{l_{1},\cdots,l_{p}=m}{\sum}}
\Big[\overset{p}{\underset{i=1}{\prod}}a_{j_{i}} \frac{C_{l_{i}}}{\sqrt{l_{i}!}}
\Big]
\\\notag
&
\times\int_{\mathbb{R}^{np}}\Big\{\overset{p}{\underset{i=1}{\prod}}
\varepsilon^{\frac{n}{\alpha}}G(\varepsilon t_{j_{i}},\varepsilon^{\frac{1}{\alpha}}\mathbf{x}_{j_{i}}-\varepsilon^{\frac{1}{\alpha}}\mathbf{y}_{i})\Big\}
\Big[\mathbb{E}\overset{p}{\underset{i=1}{\prod}}H_{l_{i}}(\zeta(\varepsilon^{-\chi}\mathbf{y}_{i}))\Big]
d\mathbf{y}_{1}\cdots d\mathbf{y}_{p}.
\end{align}

To analyze $\mathbb{E}(\xi_{\varepsilon,\leq N_{0}})^{p}$,
$p=2\nu$(the odd $p=2\nu+1$ is  unnecessary, since all the
involved random variables are centered),  we employ  the diagram
method (see, \cite{BreuerMajor} or  \cite[p.72]{Ivanov and
Leonenko}). A graph $\Gamma$ with $l_{1}+\cdots+l_{p}$ vertices is
called a (complete) diagram
of order ($l_{1},\cdots,l_{p}$) if:\\
(a) the set of vertices $V$ of the graph $\Gamma$ is of the form $V=\overset{p}{\underset{j=1}{\cup}}W_{j}$,
where $W_{j}=\{(j,l):1\leq l\leq l_{j}\}$ is the $j$-th level of the graph $\Gamma$, $1\leq j\leq p$;
\\
(b) each vertex is of degree 1; that is, each vertex is just an endpoint of an edge.
\\
(c) if $((j_{1},l_{1}),\ (j_{2},l_{2}))\in \Gamma$ then $j_{1}\neq j_{2}$; that is, the edges of the graph $\Gamma$
may connect only different levels.

Let $\mathrm{T}=\mathrm{T}(l_{1},\cdots,l_{p})$ be a set of
(complete) diagrams $\Gamma$'s of order ($l_{1},\cdots,l_{p}$).
Denote by $E(\Gamma)$  the set of edges of the graph $\Gamma\in
\mathrm{T}$. For the edge $e=((j_{1},l^{'}_{1}),\
(j_{2},l^{'}_{2}))\in E(\Gamma)$, $j_{1}<j_{2}$, $1\leq
l^{'}_{1}\leq l_{1}$, $1\leq l^{'}_{2}\leq l_{2}$, we set
$d_{1}(e)=j_{1}$, $d_{2}(e)=j_{2}$, to denote the location of the
edge $e$ in $\Gamma$. We call a diagram $\Gamma$ to be {\it
regular} if its levels can be split into pairs in such a manner
that no edge connects the levels belonging to different pairs.
Denote by $\mathrm{T}^{*}=\mathrm{T}^{*}(l_{1},\cdots,l_{p})$ the
set of all regular diagrams in $\mathrm{T}$. Therefore, if
$\Gamma\in \mathrm{T}^{*}$ is a regular diagram then it can be
divided into $p/2$ sub-diagrams (denoted by
$\Gamma_{1},\cdots,\Gamma_{p/2}$), which can not be separated
again; in this case, we naturally define $d_{1}(\Gamma_{i})\equiv
d_{1}(e)$ and $d_{2}(\Gamma_{i})\equiv d_{2}(e)$ for any $e\in
E(\Gamma_{i}),\ i=1,\ldots,\nu=p/2.$. We denote $\sharp E(\Gamma)$
(resp. $\sharp E(\Gamma_{j})$) the number of edges belonging to
the specific diagram $\Gamma$ (resp. the sub-diagram
$\Gamma_{j}$).

Based on the notations above and
let
\begin{equation*}
D_{p}=\{(J,L):J=(j_{1},\cdots,j_{p}),1\leq j_{i}\leq M,\
L=(l_{1},\cdots, l_{p}),m\leq l_{i}\leq N_{0},i=1,\cdots,p\},
\end{equation*}
(\ref{proofthm:weaksmall5}) can be rewritten as
\begin{align}\label{proofthm:weaksmall6}
\mathbb{E}(\xi_{\varepsilon,\leq N_{0}})^{p}=
\underset{(J,L)\in D_{p}}{\sum}K(J,L)\underset{\Gamma\in \mathrm{T}^{*}}{\sum}F_{\Gamma}(J,L,\varepsilon)
+\underset{(J,L)\in D_{p}}{\sum}K(J,L)\underset{\Gamma\in \mathrm{T}\backslash \mathrm{T}^{*}}{\sum}F_{\Gamma}(J,L,\varepsilon),
\end{align}
where
\begin{align}\label{proofthm:weakK(J,L)}
&K(J,L)=\overset{p}{\underset{i=1}{\prod}}a_{j_{i}} \frac{C_{l_{i}}}{\sqrt{l_{i}!}}
\\\notag
&F_{\Gamma}(J,L,\varepsilon)=\varepsilon^{-\frac{pn\chi}{2}}
\underset{\mathbb{R}^{np}}{\int}\hspace{-0.1cm}\Big\{\overset{p}{\underset{i=1}{\prod}}\varepsilon^{\frac{n}{\alpha}}G(\varepsilon t_{j_{i}},\varepsilon^{\frac{1}{\alpha}}(\mathbf{x}_{j_{i}}-\mathbf{y}_{i}))\Big\}
\hspace{-0.1cm}
\Big[
\hspace{-0.1cm}\underset{e\in E(\Gamma)}{\prod}\hspace{-0.1cm}R(\varepsilon^{-\chi}(\mathbf{y}_{d_{1}(e)}-\mathbf{y}_{d_{2}(e)}))
\Big]d\mathbf{y}_{1}\cdots d\mathbf{y}_{p}.
\end{align}
Next, we want to verify two things:
\begin{align*}
\left\{
\begin{array}{l}
(1)\ \underset{\varepsilon\rightarrow 0}{\textup{lim}}
\underset{(J,L)\in D_{p}}{\sum}K(J,L)\underset{\Gamma\in \mathrm{T}^{*}}{\sum}F_{\Gamma}(J,L,\varepsilon)
=(p-1)!!
\Big\{\mathbb{E}\Big[\big(\overset{M}{\underset{j=1}{\sum}}a_{j}V_{m,N_{0}}(t_{j},\mathbf{x}_{j})\big)^{2}\Big]\Big\}^{p/2},
\\
(2)\ \underset{\varepsilon\rightarrow 0}{\textup{lim}}\underset{(J,L)\in D_{p}}{\sum}K(J,L)\underset{\Gamma\in \mathrm{T}\backslash \mathrm{T}^{*}}{\sum}F_{\Gamma}(J,L,\varepsilon)
=0.
\end{array}\right.
\end{align*}
{\it Proof of (1):} As argued above, each $\Gamma\in\mathrm{T}^*$,
the the case $p=2\nu,\ \nu\in \mathbb{N}$, has a unique
decomposition into  sub-diagrams
$\Gamma=(\Gamma_{1},\cdots,\Gamma_{\nu})$, for which each one
cannot be further decomposed. Accordingly, we can rewrite
$F_{\Gamma}(J,L,\varepsilon)$ as the following $\nu=p/2$ products,
\begin{align}\notag
&F_{\Gamma}(J,L,\varepsilon)
\\\notag
=&
\varepsilon^{-\frac{pn\chi}{2}}
\overset{\nu}{\underset{i=1}{\prod}}\underset{\mathbb{R}^{2n}}{\int}\hspace{-0.1cm}
\varepsilon^{\frac{n}{\alpha}}
G(\varepsilon t_{d_{1}(\Gamma_{i})},\varepsilon^{\frac{1}{\alpha}}(\mathbf{x}_{d_{1}(\Gamma_{i})}-\mathbf{y}))
\varepsilon^{\frac{n}{\alpha}}
G(\varepsilon t_{d_{2}(\Gamma_{i})},\varepsilon^{\frac{1}{\alpha}}(\mathbf{x}_{d_{2}(\Gamma_{i})}-\mathbf{y}^{'}))
R^{\sharp E(\Gamma_{i})}(\varepsilon^{-\chi}(\mathbf{y}-\mathbf{y}^{'}))
\\\label{proofthm:weaksmall7}
=&
\varepsilon^{-\frac{pn\chi}{2}}
\overset{\nu}{\underset{i=1}{\prod}}\int_{\mathbb{R}^{n}}\hspace{-0.1cm}
\varepsilon^{\frac{n}{\alpha}}
G(\varepsilon (t_{d_{1}(\Gamma_{i})}+t_{d_{2}(\Gamma_{i})}),\varepsilon^{\frac{1}{\alpha}}(\mathbf{x}_{d_{1}(\Gamma_{i})}-\mathbf{x}_{d_{2}(\Gamma_{i})}-\mathbf{z}))
R^{\sharp E(\Gamma_{i})}(\varepsilon^{-\chi}\mathbf{z})
d\mathbf{z}.
\end{align}
We note that
\begin{align}\label{degeneratedMarkov1}
R^{\sharp E(\Gamma_{i})}(\varepsilon^{-\chi}\mathbf{z})
=\varepsilon^{n\chi}\int_{\mathbb{R}^{n}}e^{i<\mathbf{z},\lambda>}f^{*\sharp E(\Gamma_{i})}(\varepsilon^{\chi}\lambda)d\lambda,
\ i=1,\cdots,\nu,
\end{align}
and $\sharp E(\Gamma_{i})>n/\kappa$ ( since $\kappa> n/m$ in the
Condition B). By the Fourier transform of $G$,
\begin{align}\notag
&\int_{\mathbb{R}^{n}}e^{i<\mathbf{z},\lambda>}
\varepsilon^{\frac{n}{\alpha}}G(\varepsilon
(t_{d_{1}(\Gamma_{i})}+t_{d_{2}(\Gamma_{i})}),\varepsilon^{\frac{1}{\alpha}}(\mathbf{x}_{d_{1}(\Gamma_{i})}-\mathbf{x}_{d_{2}(\Gamma_{i})}-\mathbf{z}))
d\mathbf{z}
\\\label{modifiedgreen}
=&e^{i<\lambda,\mathbf{x}_{d_{1}(\Gamma_{i})}-\mathbf{x}_{d_{2}(\Gamma_{i})}>}
\textup{exp}\{\varepsilon(t_{d_{1}(\Gamma_{i})}+t_{d_{2}(\Gamma_{i})})
[\mathfrak{m}-(\mathfrak{m}^{\frac{2}{\alpha}}+|\varepsilon^{-\frac{1}{\alpha}}\lambda|^{2})^{\frac{\alpha}{2}}]\},
\end{align}
applying the the small-scale of $G$   illustrated in
(\ref{rescaledgreenlimit}), we have
\begin{align}\label{proofthm:weaksmall8}
\underset{\varepsilon\rightarrow 0}{\textup{lim}}F_{\Gamma}(J,L,\varepsilon)
\hspace{-0.05cm}=\hspace{-0.1cm}
\overset{\nu}{\underset{i=1}{\prod}}f^{*\sharp E(\Gamma_{i})}(\mathbf{0})
\hspace{-0.15cm}\int\hspace{-0.1cm}
e^{i<\lambda,\mathbf{x}_{d_{1}(\Gamma_{i})}-\mathbf{x}_{d_{2}(\Gamma_{i})}>}
\textup{exp}\big\{\hspace{-0.15cm}-(t_{d_{1}(\Gamma_{i})}+t_{d_{2}(\Gamma_{i})})|\lambda|^{\alpha}
\hspace{-0.05cm}\big\}d\lambda.
\end{align}
Meanwhile, if the  $L=\{l_{1},\cdots,l_{2\nu}\}$ in the defining
(\ref{proofthm:weakK(J,L)}) of $K(J,L)$ corresponds to a  regular
diagram $\Gamma$ in $\mathrm{T}(l_{1},\cdots,l_{2\nu})$, then
\begin{equation}\label{proofthm:regularcoeff}
K(J,L)=\overset{\nu}{\underset{i=1}{\prod}}
a_{d_{1}(\Gamma_{i})} a_{d_{2}(\Gamma_{i})}
\frac{C^{2}_{\sharp E(\Gamma_{i})}}{\sharp E(\Gamma_{i})!}.
\end{equation}
Therefore, by (\ref{proofthm:weaksmall8}) and (\ref{proofthm:regularcoeff}),
\begin{align}\notag
&\underset{\varepsilon\rightarrow 0}{\textup{lim}}
\underset{(J,L)\in D_{2\nu}}{\sum}K(J,L)\underset{\Gamma\in \mathrm{T}^{*}}{\sum}F_{\Gamma}(J,L,\varepsilon)
\\\notag
=&
\underset{(J,L)\in D_{2\nu}}{\sum}
\ \underset{\Gamma\in \mathrm{T}^{*}}{\sum}
\Big[
\overset{\nu}{\underset{i=1}{\prod}}
a_{d_{1}(\Gamma_{i})} a_{d_{2}(\Gamma_{i})}
\int
e^{i<\lambda,\mathbf{x}_{d_{1}(\Gamma_{i})}-\mathbf{x}_{d_{2}(\Gamma_{i})}>}
\textup{exp}\big\{\hspace{-0.15cm}-(t_{d_{1}(\Gamma_{i})}+t_{d_{2}(\Gamma_{i})})|\lambda|^{\alpha}
\big\}d\lambda
\Big]
\\\label{proofthm:weaksmall9}&\times\Big[
\overset{\nu}{\underset{i=1}{\prod}}
f^{*\sharp E(\Gamma_{i})}(\mathbf{0})
\frac{C^{2}_{\sharp E(\Gamma_{i})}}{\sharp E(\Gamma_{i})!}
\Big].
\end{align}
We note that all components in the first bracket in
$(\ref{proofthm:weaksmall9})$ are independent to the index set $L$
and the summation $\underset{\Gamma\in \mathrm{T}^{*}}{\sum}$
depends only on $\underset{L}{\sum}$, therefore
\begin{align}\notag
&\underset{\varepsilon\rightarrow 0}{\textup{lim}}
\underset{(J,L)\in D_{2\nu}}{\sum}K(J,L)\underset{\Gamma\in \mathrm{T}^{*}}{\sum}F_{\Gamma}(J,L,\varepsilon)
\\\notag
=&
\underset{L}{\sum}
\underset{\Gamma\in \mathrm{T}^{*}}{\sum}
\underset{J}{\sum}\Big[
\overset{\nu}{\underset{i=1}{\prod}}
a_{d_{1}(\Gamma_{i})} a_{d_{2}(\Gamma_{i})}
\int
e^{i<\lambda,\mathbf{x}_{d_{1}(\Gamma_{i})}-\mathbf{x}_{d_{2}(\Gamma_{i})}>}
\textup{exp}\big\{\hspace{-0.15cm}-(t_{d_{1}(\Gamma_{i})}+t_{d_{2}(\Gamma_{i})})|\lambda|^{\alpha}
\big\}d\lambda
\Big]
\\\notag&\times\Big[
\overset{\nu}{\underset{i=1}{\prod}}
f^{*\sharp E(\Gamma_{i})}(\mathbf{0})
\frac{C^{2}_{\sharp E(\Gamma_{i})}}{\sharp E(\Gamma_{i})!}
\Big]
\\\label{proofthm:regualr1}
=&
\Big[
\overset{M}{\underset{j,j^{'}=1}{\sum}}a_{j} a_{j^{'}}
\int
e^{i<\lambda,\mathbf{x}_{j}-\mathbf{x}_{j^{'}}>}
\textup{exp}\big\{\hspace{-0.15cm}-(t_{j}+t_{j^{'}})|\lambda|^{\alpha}
\big\}d\lambda
\Big]^{\nu}
\underset{L}{\sum}
\underset{\Gamma\in \mathrm{T}^{*}}{\sum}
\Big[
\overset{\nu}{\underset{i=1}{\prod}}
f^{*\sharp E(\Gamma_{i})}(\mathbf{0})
\frac{C^{2}_{\sharp E(\Gamma_{i})}}{\sharp E(\Gamma_{i})!}
\Big].
\end{align}
To handle the summation in the above, we note that $
\overset{\nu}{\underset{i=1}{\prod}} f^{*\sharp
E(\Gamma_{i})}(\mathbf{0}) \frac{C^{2}_{\sharp
E(\Gamma_{i})}}{\sharp E(\Gamma_{i})!} $ only depends on $\{\sharp
E(\Gamma_{i}),i=1,\cdots,\nu\}$, not on the structures of
sub-diagrams $\Gamma_{i},\ i=1,\cdots,\nu$; thus we may rewrite
the above summation based on the following observation. Let $s$ be
the number of different integers $r_{1},\ldots,r_{s}$ in
$\{l_{1},\cdots,l_{2\nu}\}$ with $m \leq r_{1}<\ldots <r_{s}\leq
N_{0}$. A  regular diagram requires $1\leq s\leq\nu$, which also
implies that the set $\{l_{1},\cdots,l_{2\nu}\}$ can be split into
$s$ subsets $Q_{1},\ldots,Q_{s}$ and all elements within $Q_{i}$
have the common value $r_{i},\ i=1,\ldots,s.$ For the number of
{\it pairs} within each subset $Q_{i}$, we denote it by $q_{i}$,
which satisfies $q_{i}\geq 1$, $i=1,\ldots,s$, and
$q_{1}+\cdots+q_{s}=\nu$. Thus, the above summation is
\begin{equation*}
\underset{1\leq s\leq \nu}{\sum} (s!) \ {\underset{m\leq
r_{1}<\cdots<r_{s}=N_{0}}{\sum}} \
\underset{q_{1}+\cdots+q_{s}=\nu}{\sum}\frac{(2\nu)!}{(2q_{1})!\cdots
(2q_{s})!}[\cdots].
\end{equation*}
However, for any $(s;r_{1},\ldots,r_{s};q_{1},\ldots,q_{s})$ in
the above sum,  it corresponds  $\frac{(2q_{1})!\cdots
(2q_{s})!}{2^{\nu}q_{1}!\cdots q_{s}!}(r_{1}!)^{q_{1}}\cdots
(r_{s}!)^{q_{s}}$ different regular diagrams. Therefore,
\begin{align}\notag
&\underset{L}{\sum}
\underset{\Gamma\in \mathrm{T}^{*}}{\sum}
\Big[
\overset{\nu}{\underset{i=1}{\prod}}
f^{*\sharp E(\Gamma_{i})}(\mathbf{0})
\frac{C^{2}_{\sharp E(\Gamma_{i})}}{\sharp E(\Gamma_{i})!}
\Big]
\\\notag
=&\underset{1\leq s\leq \nu}{\sum} (s!) \ {\underset{m\leq
r_{1}<\cdots<r_{s}=N_{0}}{\sum}} \
\underset{q_{1}+\cdots+q_{s}=\nu}{\sum}
\frac{(2\nu)!}{2^{\nu}q_{1}!\cdots q_{s}!}(r_{1}!)^{q_{1}}\cdots
(r_{s}!)^{q_{s}} \Big[ \overset{s}{\underset{i=1}{\prod}} \big(
f^{*r_{i}}(\mathbf{0}) \frac{C^{2}_{r_{i}}}{r_{i}!} \big)^{q_{i}}
\Big]
\\\notag
=&(2\nu-1)!! \underset{1\leq s\leq \nu}{\sum} (s!) \
{\underset{m\leq r_{1}<\cdots<r_{s}=N_{0}}{\sum}} \
\underset{q_{1}+\cdots+q_{s}=\nu}{\sum} \frac{\nu!}{q_{1}!\cdots
q_{s}!} \Big[ \overset{s}{\underset{i=1}{\prod}} \big(
f^{*r_{i}}(\mathbf{0}) C^{2}_{r_{i}} \big)^{q_{i}} \Big]
\\\label{proofthm:weaksmalla}
=&
(2\nu-1)!!
\Big[
\overset{N_{0}}{\underset{r=m}{\sum}}
f^{*r}(\mathbf{0})
C^{2}_{r}
\Big]^{\nu}.
\end{align}
Substituting (\ref{proofthm:weaksmalla}) into (\ref{proofthm:regualr1})
and recalling
$\sigma_{m,N_{0}}=(\overset{N_{0}}{\underset{r=m}{\sum}}
f^{*r}(\mathbf{0})
C^{2}_{r})^{\frac{1}{2}}$,
we get
\begin{align}\notag
&\underset{\varepsilon\rightarrow 0}{\textup{lim}}
\underset{(J,L)\in D_{2\nu}}{\sum}K(J,L)\underset{\Gamma\in \mathrm{T}^{*}}{\sum}F_{\Gamma}(J,L,\varepsilon)
\\\notag
=&
(2\nu-1)!!
\Big[
\overset{M}{\underset{j,j^{'}=1}{\sum}}a_{j} a_{j^{'}}
\int_{\mathbb{R}^{n}}
e^{i<\lambda,\mathbf{x}_{j}-\mathbf{x}_{j^{'}}>}
\textup{exp}\big\{\hspace{-0.15cm}-(t_{j}+t_{j^{'}})|\lambda|^{\alpha}
\big\}d\lambda
\Big]^{\nu}
\Big[
\overset{N_{0}}{\underset{r=m}{\sum}}
f^{*r}(\mathbf{0})
C^{2}_{r}
\Big]^{\nu}
\\\label{proofthm:regualr2}
=&
(2\nu-1)!!
\Big[
\mathbb{E}
\Big(
\overset{M}{\underset{j=1}{\sum}}a_{j}
\int_{\mathbb{R}^{n}}
e^{i<\lambda,\mathbf{x}_{j}>}
\sigma_{m,N_{0}}e^{-t_{j}|\lambda|^{\alpha}}W(d\lambda)
\Big)^{2}
\Big]^{\nu}.
\end{align}
{\it Proof of (2)}:
$
\underset{\varepsilon\rightarrow 0}{\textup{lim}}\underset{(J,L)\in D_{p}}{\sum}K(J,L)
\underset{\Gamma\in \mathrm{T}\backslash \mathrm{T}^{*}}{\sum}F_{\Gamma}(J,L,\varepsilon)
=0.$\\
By (\ref{Markovmethod2}), the number of elements in the summation
 of $\underset{(J,L)\in D_{p}}{\sum}$ is finite, thus it  suffices to
show that $\underset{\varepsilon\rightarrow 0}{\textup{lim}}
F_{\Gamma}(J,L,\varepsilon)=0$, i.e.,
\begin{align}\label{proofthm:nonregular5}
\varepsilon^{-\frac{pn\chi}{2}}\hspace{-0.2cm}
\underset{\mathbb{R}^{np}}{\int}\hspace{-0.1cm}\Big\{\overset{p}{\underset{i=1}{\prod}}\varepsilon^{\frac{n}{\alpha}}G(\varepsilon t_{j_{i}},\varepsilon^{\frac{1}{\alpha}}(\mathbf{x}_{j_{i}}-\mathbf{y}_{i}))\Big\}
\hspace{-0.1cm}
\Big[
\underset{e\in E(\Gamma)}{\prod}\hspace{-0.2cm}R(\varepsilon^{-\chi}(\mathbf{y}_{d_{1}(e)}-\mathbf{y}_{d_{2}(e)}))
\Big]d\mathbf{y}_{1}\cdots d\mathbf{y}_{p}\rightarrow0
\end{align}
for each $\Gamma\in \mathrm{T}(l_{1},\cdots,l_{p})\backslash
\mathrm{T}^{*}$. With loss of generality, we may just prove
(\ref{proofthm:nonregular5}) for  $t_{j_{i}}=1$ and
$\mathbf{x}_{j_{i}}= \mathbf{0},\ i=1,\ldots,p$, and also just
consider the case $ l_{1}\leq l_{2}\leq\ldots\leq l_{p}.$ Let
$$
A_{j,j^{'}}:=\big\{e\in E(\Gamma)\mid d_{1}(e)=j,\
d_{2}(e)=j^{'}\big\},\ B(i):=\underset{j^{'}>i}{\cup}
A_{i,j^{'}},\ 1\leq i, j<j^{'}\leq p.
$$
We observe that the number $\sharp B(i)$ of $B(i)$ must be $\le
l_i$, and a non-regular diagram $\Gamma$ must contain and an
non-empty $B(i)$ with $\sharp B(i)< l_i$; moreover, it has
(\cite[(2.20)]{BreuerMajor})
\begin{equation}
\overset{p}{\underset{i=1}{\sum}} \frac{\sharp B(i)}{l_{i}}\geq
\frac{p}{2}.
\end{equation}
\begin{align}\notag
&F_{\Gamma}(J,L,\varepsilon)
\\\notag
=&\varepsilon^{-\frac{pn\chi}{2}}
\underset{\mathbb{R}^{np}}{\int}\hspace{-0.1cm}
\Big\{\overset{p}{\underset{i=1}{\prod}}\varepsilon^{\frac{n}{\alpha}}G(\varepsilon,\varepsilon^{\frac{1}{\alpha}}\mathbf{y}_{i})\Big\}
\hspace{-0.1cm}
\Big[
\underset{i; B(i)\neq \phi}{\prod}
\ \underset{e\in B(i)}{\prod}
R(\varepsilon^{-\chi}(\mathbf{y}_{i}-\mathbf{y}_{d_{2}(e)}))
\Big]d\mathbf{y}_{1}\cdots d\mathbf{y}_{p}
\\\notag
\leq&
\varepsilon^{-\frac{pn\chi}{2}}
\underset{\mathbb{R}^{np}}{\int}\hspace{-0.1cm}
\Big\{\overset{p}{\underset{i=1}{\prod}}\varepsilon^{\frac{n}{\alpha}}G(\varepsilon,\varepsilon^{\frac{1}{\alpha}}\mathbf{y}_{i})\Big\}
\hspace{-0.1cm}
\Big[
\underset{i; B(i)\neq \phi}{\prod}
\ \underset{e\in B(i)}{\sum}
\frac{1}{\sharp B(i)}
R^{\sharp B(i)}(\varepsilon^{-\chi}(\mathbf{y}_{i}-\mathbf{y}_{d_{2}(e)}))
\Big]d\mathbf{y}_{1}\cdots d\mathbf{y}_{p}
\\\label{proofthm:nonregular1}
\leq&c
\varepsilon^{-\frac{pn\chi}{2}}
\underset{\mathbb{R}^{np}}{\int}\hspace{-0.1cm}
\Big\{\overset{p}{\underset{i=1}{\prod}}\varepsilon^{\frac{n}{\alpha}}G(\varepsilon,\varepsilon^{\frac{1}{\alpha}}\mathbf{y}_{i})\Big\}
\hspace{-0.1cm}
\Big[
\underset{i; B(i)\neq \phi}{\prod}
\ \underset{j; A_{i,j}\neq\phi}{\sum}
\frac{1}{\sharp B(i)}
R^{\sharp B(i)}(\varepsilon^{-\chi}(\mathbf{y}_{i}-\mathbf{y}_{j}))
\Big]d\mathbf{y}_{1}\cdots d\mathbf{y}_{p},
\end{align}
where
$c=\underset{i;B(i)\neq \phi}{\prod}\ \underset{j; A_{i,j}\neq\phi}{\sum}\sharp A_{i,j}/\sharp B(i)$.
\\
For any $i\in \{1,\ldots,p-1\}$ with $B(i)\neq \phi$, let  $j(i)$
be any term in $\{j^{'};A_{i,j^{'}}\neq \phi\}$. To prove
$(\ref{proofthm:nonregular1})\rightarrow 0$, by the spectral
representation, it suffices to show that
\begin{align}\label{proofthm:nonregular2}
\varepsilon^{-\frac{pn\chi}{2}}\hspace{-0.12cm}
\underset{\mathbb{R}^{np}}{\int}\hspace{-0.11cm}
\Big\{\overset{p}{\underset{i=1}{\prod}}\hspace{-0.06cm}\varepsilon^{\frac{n}{\alpha}}G(\varepsilon,\varepsilon^{\frac{1}{\alpha}}\mathbf{y}_{i})\Big\}
\hspace{-0.11cm}
\Big[
\underset{i; B(i)\neq \phi}{\prod}
\int\hspace{-0.1cm}e^{i<\mathbf{y}_{i}-\mathbf{y}_{j(i)},\lambda_{i,j(i)}>}
f^{*\sharp B(i)}(\varepsilon^{\chi}\lambda_{i,j(i)})\varepsilon^{n\chi}d\lambda_{i,j(i)}
\Big]d\mathbf{y}_{1}\cdots \hspace{-0.03cm}d\mathbf{y}_{p}
\end{align}
converges to zero when $\varepsilon\rightarrow 0$.\\
Applying Lemma 1 to $k=\sharp B(i)$, the number of $B(i)$,  we see
that,
\begin{align}\label{proofthm:nonregular3}
f^{*\sharp B(i)}(\lambda) \leq \left\{\begin{array}{lr} o(1),\ &
\textup{if}\  \sharp B(i)=l_{i},
\\
o(|\lambda|^{n(\frac{\sharp B(i)}{l_{i}}-1)}),\ & \textup{if}\
1<\sharp B(i)<l_{i},
\end{array}
\right.\ \textup{when}\ |\lambda|\rightarrow 0.
\end{align}
For example, in the case $(a)$ of Lemma 1, we can write it as
follows:
\begin{equation}\label{productionofLemma1}
f^{*\sharp B(i)}(\lambda)=C_{\sharp B(i)}(\lambda)|\lambda|^{\sharp B(i)\frac{n}{l_{i}}-n},\
C_{\sharp B(i)}(\lambda)=B_{\sharp B(i)}(\lambda)|\lambda|^{\sharp B(i)(\kappa-\frac{n}{l_{i}})},
\end{equation}
where $\underset{|\lambda|\rightarrow0}{\textup{lim}}C_{\sharp B(i)}(\lambda)=0$
because $\kappa>n/m\geq n/l_{i}$.
\\
Thus,
\begin{equation}\label{proofthm:nonregular4}
(\ref{proofthm:nonregular2})\leq\varepsilon^{-\frac{pn\chi}{2}}
o(\varepsilon^{\chi n(\sum \frac{\sharp B(i)}{l_{i}})})Q_{\varepsilon},
\end{equation}
where
\begin{equation*}
Q_{\varepsilon}=\underset{\mathbb{R}^{np}}{\int}\hspace{-0.1cm}
\Big\{\overset{p}{\underset{i=1}{\prod}}\varepsilon^{\frac{n}{\alpha}}
G(\varepsilon,\varepsilon^{\frac{1}{\alpha}}\mathbf{y}_{i})\Big\}
\hspace{-0.1cm}
\Big[
\underset{i; B(i)\neq \phi}{\prod}
\int e^{i<\mathbf{y}_{i}-\mathbf{y}_{j(i)},\lambda_{i,j(i)}>}
|\lambda_{i,j(i)}|^{n(\frac{\sharp B(i)}{l_{i}}-1)}d\lambda_{i,j(i)}
\Big]d\mathbf{y}_{1}\cdots d\mathbf{y}_{p},
\end{equation*}
which is bounded in $0<\epsilon<<1$. Because, firstly, for each
$\lambda_{i,j(i)}$, by (\ref{rescaledgreenlimit})   the following
is bounded in $0<\epsilon<<1$,
\begin{align}
\int_{\mathbb{R}^{np}}
\Big\{\overset{p}{\underset{i=1}{\prod}}\varepsilon^{\frac{n}{\alpha}}
G(\varepsilon,\varepsilon^{\frac{1}{\alpha}}\mathbf{y}_{i})\Big\}
\Big\{ \underset{i; B(i)\neq
\phi}{\prod}e^{i<\mathbf{y}_{i}-\mathbf{y}_{j(i)},\lambda_{i,j(i)}>}
\Big\} d\mathbf{y}_{1}\cdots d\mathbf{y}_{p},
\end{align}
and moreover
\begin{equation*}
\underset{i; B(i)\neq \phi}{\prod}
|\lambda_{i,j(i)}|^{n(\frac{\sharp B(i)}{l_{i}}-1)}
\end{equation*}
is integrable with respect to $\underset{i; B(i)\neq
\phi}{\prod}d\lambda_{i,j(i)}$ near the origin. Finally, the
convergence of $(\ref{proofthm:nonregular4})$ to zero is followed
by the inequality cited above, i.e. $
\overset{p}{\underset{i=1}{\sum}} \frac{\sharp B(i)}{l_{i}}\geq
\frac{p}{2}. $ \qed

\bigskip

\noindent{\bf Proof of Theorem \ref{Macro Scaling Theorem}.}\\
By the solution form (\ref{meansquaresolution}) and
$\int_{\mathbb{R}^{n}}G_{\alpha,\mathfrak{m}}(t,\mathbf{x})d\mathbf{x}=1$,
\begin{align}\notag
&T^{\frac{m\kappa}{4}}
\Big\{
u(Tt,\sqrt{T}\mathbf{x};h(\zeta(\cdot)))-C_{0}
\Big\}
\\\notag
=&
T^{\frac{m\kappa}{4}}
\Big\{
\int_{\mathbb{R}^{n}}
G_{\alpha,\mathfrak{m}}(Tt,\sqrt{T}\mathbf{x}-\mathbf{y})
\Bigl[
C_{0}
+\overset{\infty}{\underset{k=m}{\sum}}
C_{k}\frac{H_{k}(\zeta(\mathbf{y}))}{\sqrt{k!}}
\Bigr]
d\mathbf{y}-C_{0}
\Big\}\\\label{proof1dedinition}
=&
\overset{\infty}{\underset{k=m}{\sum}}
T^{\frac{m\kappa}{4}}
\frac{C_{k}}{\sqrt{k!}}
\int_{\mathbb{R}^{n}}
G_{\alpha,\mathfrak{m}}(Tt,\sqrt{T}\mathbf{x}-\mathbf{y})
H_{k}(\zeta(\mathbf{y}))d\mathbf{y}=:\overset{\infty}{\underset{k=m}{\sum}}u_{k,T}(t,\mathbf{x}).
\end{align}
By the Slutsky argument (see, for example,  \cite[p. 6.]{L}),
Theorem \ref{Macro Scaling Theorem} will be proved if we can show
that
\begin{equation}
\left\{\begin{array}{lr}
\textup{(i)}\ u_{m,T}(t,\mathbf{x}) \Rightarrow U_{m}(t,\mathbf{x}),
\\
\textup{(ii)}\overset{\infty}{\underset{k=m+1}{\sum}}u_{k,T}(t,\mathbf{x}) \rightarrow 0\ \textup{in probability},
\end{array}
\right. \textup{as}\ T\rightarrow\infty.
\end{equation}
{\it Proof of (i):} Replacing the component
$H_{m}(\zeta(\mathbf{y}))$ in the expression of
$u_{m,T}(t,\mathbf{x})$ with its It$\hat{\textup{o}}$-Wiener
expansion (\ref{itoformula}), and using  the Fourier transform
$\widehat{G}_{\alpha,\mathfrak{m}}(t,\lambda)$ of
$G_{\alpha,\mathfrak{m}}(t,\mathbf{x})$  in (\ref{fractional Green
function}), we have
\begin{align}\notag
&u_{m,T}(t,\mathbf{x})
\\\notag
=&
T^{\frac{m\kappa}{4}}\frac{C_{m}}{\sqrt{m!}}
\int_{\mathbb{R}^{n}}\hspace{-0.2cm}
G_{\alpha,\mathfrak{m}}(Tt,\sqrt{T}\mathbf{x}-\mathbf{y})
\Big\{
\int^{'}_{\mathbb{R}^{n\times m}}\hspace{-0.5cm}e^{i<\mathbf{y},\lambda_{1}+\ldots+\lambda_{m}>}\prod_{\sigma=1}^{m}{\sqrt{f(\lambda_{\sigma})}}
W(d\lambda_{\sigma})
\Big\}
d\mathbf{y}
\\\label{proofThm1second}
=&T^{\frac{m\kappa}{4}}
\frac{C_{m}}{\sqrt{m!}}
\int^{'}_{\mathbb{R}^{n\times m}}\hspace{-0.3cm}
e^{i<\sqrt{T}\mathbf{x},\lambda_{1}+\cdots+\lambda_{m}>}
\widehat{G}_{\alpha,\mathfrak{m}}(Tt,\lambda_{1}+\cdots+\lambda_{m})
\prod_{\sigma=1}^{m}{\sqrt{f(\lambda_{\sigma})}}
W(d\lambda_{\sigma}).
\end{align}
By the definition about $\int^{'}_{\mathbb{R}^{n\times m}}$ in (\ref{itoformula})
and the self-similarity property  $W(T^{-\frac{1}{2}} d\lambda)\overset{d}{=}T^{-\frac{n}{4}}W(d\lambda)$,
$u_{m,T}$ has the same finite dimensional distributions ($\overset{d}{=}$) as $\widetilde{u}_{m,T}$, where
\begin{align}\notag
\widetilde{u}_{m,T}(t,\mathbf{x})=&
\frac{C_{m}}{\sqrt{m!}}T^{\frac{m(\kappa-n)}{4}}
\int^{'}_{\mathbb{R}^{n\times m}}\hspace{-0.2cm}
e^{i<\mathbf{x},\lambda_{1}+\cdots+\lambda_{m}>}
\widehat{G}_{\alpha,\mathfrak{m}}(Tt,T^{-\frac{1}{2}}(\lambda_{1}+\cdots+\lambda_{m}))
\\\label{proofThm1second2}&\times\prod_{\sigma=1}^{m}{\sqrt{f(T^{-\frac{1}{2}}\lambda_{\sigma})}}
W(d\lambda_{\sigma}).
\end{align}
From the isometry property of the multiple Wiener integrals and
the integral representation of the limiting field $U_{m}(t,\mathbf{x})$ in (\ref{limiting field(macro)}),
\begin{align}\notag
&\mathbb{E}|\widetilde{u}_{m,T}(t,\mathbf{x})-U_{m}(t,\mathbf{x})|^{2}
\\\notag
=&C_{m}^{2}\int_{\mathbb{R}^{nm}} \Big|
T^{\frac{m(\kappa-n)}{4}}
\widehat{G}_{\alpha,\mathfrak{m}}(Tt,T^{-\frac{1}{2}}(\lambda_{1}+\cdots+\lambda_{m}))
\prod_{\sigma=1}^{m}{\sqrt{f(T^{-\frac{1}{2}}\lambda_{\sigma})}}
\\\label{proofThm1second2.1} &-
B(\mathbf{0})^{\frac{m}{2}}\frac{
\textup{exp}(-t\frac{\alpha}{2}\mathfrak{m}^{1-\frac{2}{\alpha}}|\lambda_{1}+\cdots+\lambda_{m}|^{2})
} {(|\lambda_{1}| \cdots
|\lambda_{m}|)^{\frac{n-\kappa}{2}}}\Big|^{2}\prod_{\sigma=1}^{m}{d\lambda_{\sigma}}.
\end{align}
Condition C and (\ref{proofmacroexponent})
allow us to apply the dominated convergence theorem
to show that (\ref{proofThm1second2.1}) will converge to zero when $T\rightarrow\infty$.
We note that the convergence in (\ref{proofmacroexponent})
can be shown to be monotone decreasing when $T\uparrow\infty$
for each $t>0$ and $\lambda\in \mathbb{R}^{n}$.
\\
Thus, we get
\begin{equation}\label{proofThm1second2.1end}
\underset{T\rightarrow \infty}{\textup{lim}}
\mathbb{E}|\widetilde{u}_{m,T}(t,\mathbf{x})-U_{m}(t,\mathbf{x})|^{2}=0,
\end{equation}
and the claim (i) is concluded by
$u_{m,T}\overset{d}{=}\widetilde{u}_{m,T}$ and the Cramer-Wold
theorem.
\\
{\it Proof of (ii):} By the orthogonal property
(\ref{expectionhermite}), the semigroup property of
$G_{\alpha,\mathfrak{m}}(t,\mathbf{x})$, and
(\ref{powercovariancespectralrep}), we have the following
equalities

\begin{align}\notag
&\mathbb{E}\big[(\overset{\infty}{\underset{k=m+1}{\sum}} u_{k,T}(t,\mathbf{x}))^{2}\big]
\\\notag
=&
T^{\frac{m\kappa}{2}}\overset{\infty}{\underset{k=m+1}{\sum}}C_{k}^{2}
\int_{\mathbb{R}^{n}}\int_{\mathbb{R}^{n}}
G_{\alpha,\mathfrak{m}}(Tt,\sqrt{T}\mathbf{x}-\mathbf{y})
G_{\alpha,\mathfrak{m}}(Tt,\sqrt{T}\mathbf{x}-\mathbf{y}^{'})
R^{k}(\mathbf{y}-\mathbf{y}^{'})d\mathbf{y}\ d\mathbf{y}^{'}
\\\notag
=&
T^{\frac{m\kappa}{2}}\overset{\infty}{\underset{k=m+1}{\sum}}C_{k}^{2}
\int_{\mathbb{R}^{n}}
G_{\alpha,\mathfrak{m}}(2Tt,\mathbf{z})
R^{k}(\mathbf{z})d\mathbf{z}
\\\notag
=&
T^{\frac{m\kappa}{2}}\overset{\infty}{\underset{k=m+1}{\sum}}C_{k}^{2}
\int_{\mathbb{R}^{n}}
\widehat{G}_{\alpha,\mathfrak{m}}(2Tt,\lambda)
f^{*k}(\lambda)d\lambda\ \ \ \ (\textup{by}\ (\ref{powercovariancespectralrep}))
\\\label{proofmacrodegeneratepart}
=&
T^{\frac{m\kappa-n}{2}}
\Big(\overset{k^{*}}{\underset{k=m+1}{\sum}}+\overset{\infty}{\underset{k=k^{*}+1}{\sum}}\Big)
C_{k}^{2}
\int_{\mathbb{R}^{n}}\hspace{-0.15cm}
\widehat{G}_{\alpha,\mathfrak{m}}(2Tt,T^{-\frac{1}{2}}\lambda)
f^{*k}(T^{-\frac{1}{2}}\lambda)d\lambda=:(I)+(II),
\end{align}
where $k^{*}=\textup{max}\{k\in \mathbb{N}|\ k\geq m+1,\ k\kappa\leq n\}$.
\\
For the case $k^{*}\kappa<n$,
by Lemma 1 and (\ref{proofmacroexponent}),
\begin{align}\notag
\underset{T\rightarrow \infty}{\textup{lim}}(I)
=&
\underset{T\rightarrow \infty}{\textup{lim}}
T^{\frac{m\kappa-n}{2}}
\overset{k^{*}}{\underset{k=m+1}{\sum}}
C_{k}^{2}
\int_{\mathbb{R}^{n}}
\widehat{G}_{\alpha,\mathfrak{m}}(2Tt,T^{-\frac{1}{2}}\lambda)
B_{k}(T^{-\frac{1}{2}}\lambda)|T^{-\frac{1}{2}}\lambda|^{k\kappa-n}d\lambda
\\\notag
\leq&
\underset{T\rightarrow \infty}{\textup{lim}}
\overset{k^{*}}{\underset{k=m+1}{\sum}}T^{\frac{m\kappa-k\kappa}{2}}
C_{k}^{2}\parallel\hspace{-0.1cm}B_{k}\hspace{-0.1cm}\parallel_{\infty}
\int_{\mathbb{R}^{n}}
e^{-t\frac{\alpha}{2}\mathfrak{m}^{1-\frac{2}{\alpha}}|\lambda|^{2}}
|\lambda|^{k\kappa-n}d\lambda
\\\notag
\leq&
\underset{T\rightarrow \infty}{\textup{lim}}
T^{-\frac{\kappa}{2}}\overset{k^{*}}{\underset{k=m+1}{\sum}}
C_{k}^{2}\parallel\hspace{-0.1cm}B_{k}\hspace{-0.1cm}\parallel_{\infty}
\int_{\mathbb{R}^{n}}
e^{-t\frac{\alpha}{2}\mathfrak{m}^{1-\frac{2}{\alpha}}|\lambda|^{2}}
|\lambda|^{k\kappa-n}d\lambda
=0.
\end{align}
For the case $k^{*}\kappa=n$,
we still have
$\underset{T\rightarrow \infty}{\textup{lim}}(I)=0$
because
\begin{equation*}
\underset{T\rightarrow \infty}{\textup{lim}}
T^{\frac{m\kappa-n}{2}}
C_{k^{*}}^{2}
\int_{\mathbb{R}^{n}}
\widehat{G}_{\alpha,\mathfrak{m}}(2Tt,T^{-\frac{1}{2}}\lambda)
B_{k^{*}}(T^{-\frac{1}{2}}\lambda)\textup{ln}(2+T^{\frac{1}{2}}|\lambda|^{-1})d\lambda=0.
\end{equation*}
On the other hand, by the assumption $m\kappa<n$ in Condition C and Lemma 1, for any $k>k^{*}+1$ we have
$
\parallel\hspace{-0.1cm} f^{*k}\hspace{-0.1cm}\parallel_{\infty}
\leq
\parallel\hspace{-0.1cm} f^{*(k^{*}+1)}\hspace{-0.1cm}\parallel_{\infty}
$, so
\begin{equation*}
\underset{T\rightarrow \infty}{\textup{lim}}(II)
\leq
\underset{T\rightarrow \infty}{\textup{lim}}
T^{\frac{m\kappa-n}{2}}
\overset{\infty}{\underset{k=k^{*}+1}{\sum}}
C_{k}^{2}\parallel\hspace{-0.1cm} f^{*(k^{*}+1)}\hspace{-0.1cm}\parallel_{\infty}
\int_{\mathbb{R}^{n}}
\widehat{G}_{\alpha,\mathfrak{m}}(2Tt,T^{-\frac{1}{2}}\lambda)
d\lambda=0.
\end{equation*}
Hence $\underset{T\rightarrow \infty}{\textup{lim}}
\mathbb{E}\big[(\overset{\infty}{\underset{k=m+1}{\sum}}
u_{k,T}(t,\mathbf{x}))^{2}\big]=0 $ and the claim (ii) is  proved
by the Markov inequality. \qed
%


\bigskip

\noindent{\bf Proof of Theorem \ref{Micro Scaling Theorem}.}
\\
The following proof is a hybrid of the proofs of Theorems 2 and 3,
we give a full presentation mainly to see how the rescaling of the
initial data is proceeded. By the Hermite expansion and the
solution form (\ref{meansquaresolution}) we can rewrite
\begin{align}\label{expansion single micro}
u^{\varepsilon}(t,\mathbf{x})=\hspace{-0.1cm}
\overset{\infty}{\underset{k=m}{\sum}}
\varepsilon^{-\frac{\chi m\kappa}{2}}\frac{C_{k}}{\sqrt{k!}}
\underset{\mathbb{R}^{n}}{\int}G_{\alpha,\mathfrak{m}}(\varepsilon
t,\mathbf{y})H_{k}(\zeta(\varepsilon^{-\frac{1}{\alpha}-\chi}(\varepsilon^{\frac{1}{\alpha}}\mathbf{x}-\mathbf{y})))d\mathbf{y}
=:\hspace{-0.1cm}
\overset{\infty}{\underset{k=m}{\sum}}I^{\varepsilon}_{k}(t,\mathbf{x}).
\end{align}
By the Slutsky argument again,  we  show that
\begin{equation}
\left\{\begin{array}{lr}
\textup{(i)}\ I^{\varepsilon}_{m}(t,\mathbf{x}) \Rightarrow V_{m}(t,\mathbf{x}),
\\
\textup{(ii)}\overset{\infty}{\underset{k=m+1}{\sum}}I^{\varepsilon}_{k}(t,\mathbf{x}) \rightarrow 0\ \textup{in probability},
\end{array}
\right. \textup{as}\ \varepsilon\rightarrow 0.
\end{equation}
{\it Proof of (i):}
By substituting the It$\hat{\textup{o}}$-Wiener expansion (\ref{itoformula})
for the random field $H_{m}(\zeta(\cdot))$ into $I^{\varepsilon}_{m}(t,\mathbf{x})$
and exchanging the order of integration
\begin{align}\notag
&I^{\varepsilon}_{m}(t,\mathbf{x}) \\\notag
=&
\frac{C_{m}}{\sqrt{m!}}\varepsilon^{-\frac{\chi m\kappa}{2}}\underset{\mathbb{R}^{n}}{\int}G_{\alpha,\mathfrak{m}}(\varepsilon
t,\mathbf{y})H_{m}(\zeta(\varepsilon^{-\frac{1}{\alpha}-\chi}(\varepsilon^{\frac{1}{\alpha}}\mathbf{x}-\mathbf{y})))d\mathbf{y}
\\\notag
=& \frac{C_{m}}{\sqrt{m!}}\varepsilon^{-\frac{\chi m\kappa}{2}}\underset{\mathbb{R}^{n}}{\int}G_{\alpha,\mathfrak{m}}(\varepsilon t,\mathbf{y})
\underset{\mathbb{R}^{n\times m}}{\int^{'}}
e^{i<\varepsilon^{-\frac{1}{\alpha}-\chi}(\varepsilon^{\frac{1}{\alpha}}\mathbf{x}-\mathbf{y}),\lambda_{1}+\cdots+\lambda_{m}>}
\overset{m}{\underset{{\sigma=1}}{\prod}}{\sqrt{f(\lambda_{\sigma})}}
W(d\lambda_{\sigma})d\mathbf{y}
\\\notag
=&\frac{C_{m}}{\sqrt{m!}}\varepsilon^{-\frac{\chi m\kappa}{2}}\underset{\mathbb{R}^{n\times m}}{\int^{'}}
e^{i<\varepsilon^{-\chi}\mathbf{x},\lambda_{1}+\cdots+\lambda_{m}>}
\widehat{G}_{\alpha,\mathfrak{m}}
(\varepsilon t,\varepsilon^{-\frac{1}{\alpha}-\chi}(\lambda_{1}+\cdots+\lambda_{m}))
\overset{m}{\underset{{\sigma=1}}{\prod}}{\sqrt{f(\lambda_{\sigma})}}
W(d\lambda_{\sigma})d\mathbf{y}
\\\notag
\overset{d}{=}& \frac{C_{m}}{\sqrt{m!}}\varepsilon^{\frac{\chi m(n-\kappa)}{2}}
\underset{{\mathbb{R}^{n\times m}}}{\int^{'}}
e^{i<\mathbf{x},\lambda^{'}_{1}+\cdots+\lambda^{'}_{m}>}
\widehat{G}_{\alpha,\mathfrak{m}}
(\varepsilon t,\varepsilon^{-\frac{1}{\alpha}}(\lambda^{'}_{1}+\cdots+\lambda^{'}_{m}))
\prod_{\sigma=1}^{m}{\sqrt{f(\varepsilon^{\chi}\lambda^{'}_{\sigma})}}
W(d\lambda^{'}_{\sigma})
\\\label{proof micro single add1}=:&\widetilde{I}^{\varepsilon}_{m}(t,\mathbf{x}),
\end{align}
where we have used the self-similarity property
$W(\varepsilon^{\chi} d\lambda)\overset{d}{=}\varepsilon^{\frac{n\chi}{2}}W(d\lambda)$
in the last equality.

Now, applying the isometry property of the multiple Wiener integrals to
the difference of $\widetilde{I}^{\varepsilon}_{m}(t,\mathbf{x})$
and the random field $V_{m}(t,\mathbf{x})$
in (\ref{limiting field(micro)}), we have
\begin{align}\notag
\mathbb{E}|\widetilde{I}^{\varepsilon}_{m}(t,\mathbf{x}) -V_{m}(t,\mathbf{x})|^{2}
=&
C_{m}^{2}\int_{\mathbb{R}^{nm}}
\big|\varepsilon^{\frac{\chi m(n-\kappa)}{2}}
\widehat{G}_{\alpha,\mathfrak{m}}(\varepsilon t,\varepsilon^{-\frac{1}{\alpha}}(\lambda_{1}+\cdots+\lambda_{m}))
\prod_{\sigma=1}^{m}{\sqrt{f(\varepsilon^{\chi}\lambda_{\sigma})}}
\\\label{proofmicrolimitlL2}&-
B(\mathbf{0})^{\frac{m}{2}}e^{- t
|\lambda_{1}+\cdots+\lambda_{m}|^{\alpha}}(|\lambda_{1}|\cdots|\lambda_{m}|)^{\frac{\kappa-n}{2}}
\big|^{2}
\prod_{\sigma=1}^{m}{d\lambda_{\sigma}}\rightarrow0
\end{align}
when $\varepsilon\rightarrow0$, by Condition C and
(\ref{rescaledgreenlimit}).
\\
By the Markov inequality, (\ref{proofmicrolimitlL2}) implies
$\widetilde{I}^{\varepsilon}_{m}(t,\mathbf{x})\rightarrow
V_{m}(t,\mathbf{x})$ in probability. However, because
$I^{\varepsilon}_{m}(t,\mathbf{x})\overset{d}{=}\widetilde{I}^{\varepsilon}_{m}(t,\mathbf{x})$,
the claim (i) is concluded,  by the Cramer-Wold argument.
\\
{\it Proof of (ii):}
From (\ref{expansion single micro}), by the orthogonal property (\ref{expectionhermite})
and the semigroup property of $G_{\alpha,\mathfrak{m}}(t,\mathbf{x})$,
\begin{align}\notag
&\mathbb{E}\Big(\overset{\infty}{\underset{k= m+1}{\sum}}I^{\varepsilon}_{k}(t,\mathbf{x})\Big)^{2}
=
\overset{\infty}{\underset{k= m+1}{\sum}}\mathbb{E}(I^{\varepsilon}_{k}(t,\mathbf{x}))^{2}
\\\notag=&
\overset{\infty}{\underset{k= m+1}{\sum}}
\varepsilon^{-\chi m\kappa}C_{k}^{2}
\int_{\mathbb{R}^{n}}\int_{\mathbb{R}^{n}}
G_{\alpha,\mathfrak{m}}(\varepsilon
t,\mathbf{y})
G_{\alpha,\mathfrak{m}}(\varepsilon
t,\mathbf{y^{'}})
R^{k}(\varepsilon^{-\frac{1}{\alpha}-\chi}(\mathbf{y}-\mathbf{y}^{'}))d\mathbf{y}d\mathbf{y}^{'}
\\\notag
=&\overset{\infty}{\underset{k= m+1}{\sum}}
\varepsilon^{-\chi m\kappa}C_{\rho}^{2}
\int_{\mathbb{R}^{n}}
G_{\alpha,\mathfrak{m}}(2\varepsilon
t,\mathbf{z})
R^{k}(\varepsilon^{-\frac{1}{\alpha}-\chi}\mathbf{z})d\mathbf{z}
\\\notag
=&
\overset{\infty}{\underset{k= m+1}{\sum}}
\varepsilon^{-\chi m\kappa}C_{k}^{2}
\int_{\mathbb{R}^{n}}
\widehat{G}_{\alpha,\mathfrak{m}}(2\varepsilon
t,\varepsilon^{-\frac{1}{\alpha}-\chi}\lambda)
f^{*k}(\lambda)d\lambda
\\\notag
=&
\Big(\overset{k^{*}}{\underset{k=m+1}{\sum}}+\overset{\infty}{\underset{k=k^{*}+1}{\sum}}\Big)
\varepsilon^{\chi(n-m\kappa)}C_{k}^{2}
\int_{\mathbb{R}^{n}}
\widehat{G}_{\alpha,\mathfrak{m}}(2\varepsilon
t,\varepsilon^{-\frac{1}{\alpha}}\lambda)
f^{*k}(\varepsilon^{\chi}\lambda)d\lambda=:(I)+(II),
\end{align}
where $k^{*}=\textup{max}\{k\in \mathbb{N}|\ k\geq m+1,\ k\kappa\leq n\}$.
\\
For the case $k^{*}\kappa<n$, by Lemma 1,
\begin{align*}
\underset{\varepsilon\rightarrow 0}{\textup{lim}}(I)
=&
\underset{\varepsilon\rightarrow 0}{\textup{lim}}
\overset{k^{*}}{\underset{k=m+1}{\sum}}
\varepsilon^{\chi(n-m\kappa)}C_{k}^{2}
\int_{\mathbb{R}^{n}}
\widehat{G}_{\alpha,\mathfrak{m}}(2\varepsilon
t,\varepsilon^{-\frac{1}{\alpha}}\lambda)
B_{k}(\varepsilon^{\chi}\lambda)|\varepsilon^{\chi}\lambda|^{k\kappa-n}
d\lambda
\\
\leq&
\underset{\varepsilon\rightarrow 0}{\textup{lim}}
\overset{k^{*}}{\underset{k=m+1}{\sum}}
\varepsilon^{\chi\kappa(k-m)}C_{k}^{2}
\parallel\hspace{-0.15cm}B_{k}\hspace{-0.15cm}\parallel_{\infty}
\int_{\mathbb{R}^{n}}
e^{-2t|\lambda|^{\alpha}}
|\lambda|^{k\kappa-n}
d\lambda=0.
\end{align*}
For the case, $k^{*}\kappa=n$, we still have
$\underset{\varepsilon\rightarrow 0}{\textup{lim}}(I)=0$
because
\begin{equation*}
\underset{\varepsilon\rightarrow 0}{\textup{lim}}
\varepsilon^{\chi(n-m\kappa)}C_{k^{*}}^{2}
\int_{\mathbb{R}^{n}}
\widehat{G}_{\alpha,\mathfrak{m}}(2\varepsilon
t,\varepsilon^{-\frac{1}{\alpha}}\lambda)
B_{k^{*}}(\varepsilon^{\chi}\lambda)\textup{ln}(2+|\varepsilon^{\chi}\lambda|^{-1})
d\lambda=0.
\end{equation*}
On the other hand, by the assumption $\kappa<n/m$ in Condition C and Lemma 1, for any $k>k^{*}+1$ we have
$
\parallel\hspace{-0.1cm} f^{*k}\hspace{-0.1cm}\parallel_{\infty}
\leq
\parallel\hspace{-0.1cm} f^{*(k^{*}+1)}\hspace{-0.1cm}\parallel_{\infty}
$, so
\begin{equation*}
\underset{\varepsilon\rightarrow 0}{\textup{lim}}\ (II)
\leq
\underset{\varepsilon\rightarrow 0}{\textup{lim}}
\overset{\infty}{\underset{k=k^{*}+1}{\sum}}
\varepsilon^{\chi(n-m\kappa)}C_{k}^{2}
\parallel\hspace{-0.1cm} f^{*(k^{*}+1)}\hspace{-0.1cm}\parallel_{\infty}
\int_{\mathbb{R}^{n}}
e^{-2t|\lambda|^{\alpha}}
d\lambda=0.
\end{equation*}
Hence $\underset{\varepsilon\rightarrow 0}{\textup{lim}}\
\mathbb{E}\big[(\overset{\infty}{\underset{k=m+1}{\sum}}
I_{k}^{\varepsilon}(t,\mathbf{x}))^{2}\big]=0 $ and the claim (ii)
is  proved by the Markov inequality. \qed

\bigskip

\noindent {\bf Appendix: Proof of Lemma 1.}

The idea of following proofs comes from \cite[p.115, Theorem
3]{Smirnov} and \cite[p.160, Theorem 8.8]{Krasnoselskii}. We only
consider  their results for the density functions on  the whole
space. Suppose that   two spectral density functions $f_{1}$ and
$f_{2}$ are in the form
\begin{equation}\label{proofsingularconvolution2}
0\leq
f_{j}(\lambda)=\frac{K_{j}(\lambda)}{|\lambda|^{n-\kappa_{j}}},\
\kappa_{j}>0,\ j=1,2,
\end{equation}
where $K_{1}(\lambda)$ and $K_{2}(\lambda)$ are nonnegative functions belonging to $\mathrm{C}(\mathbb{R}^{n}\backslash\{\mathbf{0}\})$.
\\
 Let
$g(\lambda)=\int_{\mathbb{R}^{n}}f_{1}(\lambda-\eta)f_{2}(\eta)d\eta,\
\lambda\in \mathbb{R}^{n}. $ To prove Lemma 1, we  show that $g$
can be written as
\begin{equation*}
g(\lambda)=
\left\{\begin{array}{lr}
B(\lambda)|\lambda|^{\kappa_{1}+\kappa_{2}-n},\ \ & \textup{for}\ \kappa_{1}+\kappa_{2}<n,
\\
B(\lambda)\textup{ln}(2+\frac{1}{|\lambda|}),\ \ & \textup{for}\ \kappa_{1}+\kappa_{2}=n,
\\
B(\lambda)\in \mathrm{C}(\mathbb{R}^{n}),\ \ & \textup{for}\ \kappa_{1}+\kappa_{2}> n,
\end{array}
\right.
\end{equation*}
for some bounded function $B(\lambda)\in \mathrm{C}(\mathbb{R}^{n}\backslash\{\mathbf{0}\})$.
\\
{\it Case 1: $\kappa_{1}+\kappa_{2}<n$.}
For any $\lambda_{0}\neq \mathbf{0}$,
we divide $\mathbb{R}^{n}$ into four parts: $\mathbb{R}^{n}=D_{1}\cup D_{2}\cup D_{3}\cup D_{4}$,
where
\begin{equation*}
\begin{array}{l}
D_{1}=\big\{\eta\in \mathbb{R}^{n}|\ |\eta-\lambda_{0}|<|\lambda_{0}|/2\big\},
\\
D_{2}=\big\{\eta\in \mathbb{R}^{n}|\ |\eta|<|\lambda_{0}|/2\big\},
\\
D_{3}=\big\{\eta\in (D_{1}\cup D_{2})^{c}|\ |\eta-\lambda_{0}|<|\eta|\big\},
\\
D_{4}=\big\{\eta\in (D_{1}\cup D_{2})^{c}|\ |\eta-\lambda_{0}|>|\eta|\big\}.
\end{array}
\end{equation*}
Therefore,
\begin{equation*}
g(\lambda_{0})=\overset{4}{\underset{j=1}{\sum}}
\int_{D_{j}}f_{1}(\lambda_{0}-\eta)f_{2}(\eta)d\eta=:I_{1}+I_{2}+I_{3}+I_{4}.
\end{equation*}
\begin{align*}
I_{1}(\lambda_{0})\leq&
\big(\underset{\eta^{'}\in D_{1}}{\textup{sup}}f_{2}(\eta^{'})\big)
\int_{D_{j}}f_{1}(\lambda_{0}-\eta)d\eta
\\
\leq&
\big(\underset{\eta\in D_{1}}{\textup{sup}}K_{2}(\eta)\big)(\frac{|\lambda_{0}|}{2})^{\kappa_{2}-n}
\big(\underset{\eta\in D_{2}}{\textup{sup}}K_{1}(\eta)\big)c_{n}
\int_{0}^{\frac{|\lambda_{0}|}{2}}r^{\kappa_{1}-1}dr
=C|\lambda_{0}|^{\kappa_{1}+\kappa_{2}-n}
\end{align*}
where $c_{n}$ is the surface area of the unit sphere on $\mathbb{R}^{n}$
and $C$ is a constant independent to $\lambda_{0}$.
Similarly, $I_{2}(\lambda_{0})\leq C|\lambda_{0}|^{\kappa_{1}+\kappa_{2}-n}$.
\\
By the fact $I_{3}\in\mathrm{C}(\mathbb{R}^{n}\backslash\{\mathbf{0}\})\cap L^{1}(\mathbb{R}^{n})$,
we know $\underset{|\lambda_{0}|\geq1}{\textup{sup}}I_{3}(\lambda_{0})<\infty$.
So
we suffice to study the behavior of $I_{3}(\cdot)$ on the domain $\{\lambda_{0}|\  |\lambda_{0}|<1\}$.

By the requirement (\ref{proofsingularconvolution2}),
$\underset{|\eta|\rightarrow\infty}{\textup{lim}}
K_{j}(\eta)|\eta|^{\kappa_{j}}=0$; that is, for any $\varepsilon>0$,
there exists a constant $M=M(\varepsilon)>0$ such that
\begin{align}\label{decayK}
K_{j}(\eta)\leq \varepsilon|\eta|^{-\kappa_{j}}\ \ \textup{for all}\ |\eta|>M.
\end{align}
Because $|\eta-\lambda_{0}|<|\eta|$ for $\eta\in D_{3}$,
\begin{align*}
I_{3}(\lambda_{0})\leq
\big(\hspace{-0.12cm}
\int_{D_{3}\cap \{|\eta-\lambda_{0}|>M+1\}}
+\hspace{-0.1cm}
\int_{D_{3}\cap \{|\eta-\lambda_{0}|<M+1\}}
\big)
\frac{K_{1}(\lambda_{0}-\eta)K_{2}(\eta)}{|\lambda_{0}-\eta|^{2n-\kappa_{1}-\kappa_{2}}}
d\eta
=:I_{3,1}(\lambda_{0})+I_{3,2}(\lambda_{0}).
\end{align*}
By using (\ref{decayK}) and $|\eta-\lambda_{0}|<|\eta|$ again,
\begin{align}\notag
I_{3,1}(\lambda_{0})
\leq&
\varepsilon\int_{D_{3}\cap\{|\eta-\lambda_{0}|>M+1\}}
\frac{K_{1}(\lambda_{0}-\eta)|\eta|^{-\kappa_{2}}}{|\lambda_{0}-\eta|^{2n-\kappa_{1}-\kappa_{2}}}
d\eta
\leq
\varepsilon\int_{\{|\eta-\lambda_{0}|>M+1\}}
\frac{K_{1}(\lambda_{0}-\eta)}{|\lambda_{0}-\eta|^{2n-\kappa_{1}}}
d\eta
\\\label{proofsingularconvolutioncase1}\leq&
\varepsilon(M+1)^{-n}\int_{\mathbb{R}^{n}}
\frac{K_{1}(\eta)}{|\eta|^{n-\kappa_{1}}}
d\eta=\varepsilon(M+1)^{-n}
\end{align}

\begin{align}\label{proofsingularconvolutioncase1s}
I_{3,2}(\lambda_{0})\leq \parallel \hspace{-0.15cm}K_{1}\hspace{-0.15cm}\parallel_{\infty}
\parallel \hspace{-0.15cm}K_{2}\hspace{-0.15cm}\parallel_{\infty}c_{n}
\int_{\frac{\lambda_{0}}{2}}^{M+1}
\frac{r^{n-1}}{r^{2n-\kappa_{1}-\kappa_{2}}}
d\eta
<
\frac{\parallel \hspace{-0.15cm}K_{1}\hspace{-0.15cm}\parallel_{\infty}
\parallel \hspace{-0.15cm}K_{2}\hspace{-0.15cm}\parallel_{\infty}c_{n}}{n-\kappa_{1}-\kappa_{2}}
(\frac{|\lambda_{0}|}{2})^{\kappa_{1}+\kappa_{2}-n}.
\end{align}
Combining (\ref{proofsingularconvolutioncase1}) and (\ref{proofsingularconvolutioncase1s}), we get
$I_{3}(\lambda)=B(\lambda)|\lambda|^{\kappa_{1}+\kappa_{2}-n}$
for some bounded function $B$. This observation still holds for $I_{4}$.
Therefore, the proof for the case $\kappa_{1}+\kappa_{2}<n$ is finished.
\\
{\it Case 2: $\kappa_{1}+\kappa_{2}=n,\ \kappa_{1},\ \kappa_{2}>0$.}
Let $\hat{\lambda}_{0}=\lambda_{0}/|\lambda_{0}|,$
\begin{align}\notag
g(\lambda_{0})
=&
\int_{\mathbb{R}^{n}}\frac{K_{1}(\lambda_{0}-\eta)K_{2}(\eta)}
{|\lambda_{0}-\eta|^{n-\kappa_{1}}|\eta|^{n-\kappa_{2}}}d\eta
=
\big(\int_{\{|\eta|<2\}}+\int_{|\eta|>2}\big)
\frac{K_{1}(|\lambda_{0}|(\widehat{\lambda}_{0}-\eta))K_{2}(|\lambda_{0}|\eta)}
{|\widehat{\lambda}_{0}-\eta|^{n-\kappa_{1}}|\eta|^{n-\kappa_{2}}}d\eta
\\\notag=:&J_{1}(\lambda_{0})+J_{2}(\lambda_{0}).
\end{align}
\begin{equation}\label{appendixlog1}
J_{1}(\lambda_{0})\leq
\int_{\{|\eta|<2\}}
\frac{\parallel \hspace{-0.15cm}K_{1}\hspace{-0.15cm}\parallel_{\infty}
\parallel \hspace{-0.15cm}K_{2}\hspace{-0.15cm}\parallel_{\infty}
d\eta}
{|\widehat{\lambda}_{0}-\eta|^{n-\kappa_{1}}|\eta|^{n-\kappa_{2}}}
=
\int_{\{|\eta|<2\}}
\frac{\parallel \hspace{-0.15cm}K_{1}\hspace{-0.15cm}\parallel_{\infty}
\parallel \hspace{-0.15cm}K_{2}\hspace{-0.15cm}\parallel_{\infty}
d\eta}
{|\widehat{\mathbf{x}}-\eta|^{n-\kappa_{1}}|\eta|^{n-\kappa_{2}}}
<\infty,
\end{equation}
where the last equality holds for any unit vector $\widehat{\mathbf{x}}$.
\begin{align*}
J_{2}(\lambda_{0})=
\big(\hspace{-0.2cm}
\underset{{\{2<|\eta|<2(2+\frac{1}{|\lambda_{0}|})\}}}{\int}
\hspace{-0.2cm}+\hspace{-0.2cm}\underset{\{|\eta|>2(2+\frac{1}{|\lambda_{0}|})\}}{\int}
\hspace{-0.2cm}\big)
\frac{K_{1}(|\lambda_{0}|(\widehat{\lambda}_{0}-\eta))K_{2}(|\lambda_{0}|\eta)}
{|\widehat{\lambda}_{0}-\eta|^{n-\kappa_{1}}|\eta|^{n-\kappa_{2}}}d\eta
=:J_{2,1}(\lambda_{0})+J_{2,2}(\lambda_{0}).
\end{align*}
Because $|\widehat{\lambda}_{0}-\eta|\geq |\eta|-1$
\begin{align}\notag
J_{2,1}(\lambda_{0})\leq
\underset{{\{2<|\eta|<2(2+\frac{1}{|\lambda_{0}|})\}}}{\int}
\frac{\parallel \hspace{-0.15cm}K_{1}\hspace{-0.15cm}\parallel_{\infty}
\parallel \hspace{-0.15cm}K_{2}\hspace{-0.15cm}\parallel_{\infty}
d\eta}{(|\eta|-1)^{n-\kappa_{1}}|\eta_{2}|^{n-\kappa_{2}}}
=
c_{n}\int_{2}^{2(2+\frac{1}{|\lambda_{0}|})}
\frac{\parallel \hspace{-0.15cm}K_{1}\hspace{-0.15cm}\parallel_{\infty}
\parallel \hspace{-0.15cm}K_{2}\hspace{-0.15cm}\parallel_{\infty}
dr}{(r-1)^{n-\kappa_{1}}r^{1-\kappa_{2}}}
\\\label{appendixlog2}
=c_{n}
\parallel \hspace{-0.15cm}K_{1}\hspace{-0.15cm}\parallel_{\infty}
\parallel \hspace{-0.15cm}K_{2}\hspace{-0.15cm}\parallel_{\infty}
\int_{2}^{2(2+\frac{1}{|\lambda_{0}|})}
\hspace{-0.2cm}\frac{dr}{(1-\frac{1}{r})^{n-\kappa_{1}}r}
\leq
2^{n-\kappa_{1}}c_{n}
\parallel \hspace{-0.15cm}K_{1}\hspace{-0.15cm}\parallel_{\infty}
\parallel \hspace{-0.15cm}K_{2}\hspace{-0.15cm}\parallel_{\infty}
\textup{ln}(2+\frac{1}{|\lambda_{0}|}).
\end{align}
Changing from variable $\eta$ to $\frac{\tau}{|\lambda_{0}|}$
and using the inequality
$|\lambda_{0}-\tau|\geq |\tau|-|\lambda_{0}|\geq 2(1+2|\lambda_{0}|)-|\lambda_{0}|\geq 2$
for $|\tau|>2(1+|\lambda_{0}|)$,
\begin{equation*}
J_{2,2}(\lambda_{0})=\hspace{-0.15cm}
\underset{\{|\tau|>2(1+2|\lambda_{0}|)\}}{\int}
\hspace{-0.2cm}\frac{K_{1}(\lambda_{0}-\tau)K_{2}(\eta)}
{|\lambda_{0}-\tau|^{n-\kappa_{1}}|\tau|^{n-\kappa_{2}}}d\tau
\leq
\frac{\parallel\hspace{-0.15cm}K_{1}\hspace{-0.15cm}\parallel_{\infty}}{2^{n-\kappa_{1}}}
\hspace{-0.22cm}\underset{\{|\tau|>2(1+2|\lambda_{0}|)\}}{\int}\hspace{-0.2cm}
\frac{K_{2}(\tau)}
{|\tau|^{n-\kappa_{2}}}d\eta
\leq
\frac{\parallel\hspace{-0.15cm}K_{1}\hspace{-0.15cm}\parallel_{\infty}}{2^{n-\kappa_{1}}}.
\end{equation*}
The last estimation, together with (\ref{appendixlog1}) and (\ref{appendixlog2}),
implies that
there exist bounded and positive functions $\widetilde{B}(\lambda)$ and $C(\lambda)$
such that
$g(\lambda)=\widetilde{B}(\lambda)\textup{ln}(2+\frac{1}{|\lambda|})+C(\lambda)
=B(\lambda)\textup{ln}(2+\frac{1}{|\lambda|})$,
where $B(\lambda)=\widetilde{B}(\lambda)+\frac{C(\lambda)}{\textup{ln}(2+|\lambda|^{-1})}$
is also a bounded function.\\
{\it Case 3: $\kappa_{1}+\kappa_{2}>n,\ \kappa_{1},\ \kappa_{2}>0.$}
Because $\kappa_{1}+\kappa_{2}>n$ implies that
there exist $p,\ p^{'}\in(1,\infty)$ such that
$p(n-\kappa_{1}),\ p^{'}(n-\kappa_{2})<n$ and $\frac{1}{p}+\frac{1}{p^{'}}=1$.
For any $\lambda\in \mathbb{R}^{n}$,
by H$\ddot{\textup{o}}$lder's inequality,
$g(\lambda)\leq \parallel\hspace{-0.12cm}f_{1}\hspace{-0.12cm}\parallel_{p}
\parallel\hspace{-0.12cm}f_{2}\hspace{-0.12cm}\parallel_{p^{'}}.$
Meanwhile,
it also implies
the continuity of $g$ as follows:
\begin{equation*}
|g(\lambda)-g(\lambda_{0})|=|\int_{\mathbb{R}^{n}}\hspace{-0.12cm}\big( f_{1}(\lambda-\eta)-f_{1}(\lambda_{0}-\eta)\big)f_{2}(\eta)d\eta|
\leq
\parallel\hspace{-0.13cm}f_{1}(\lambda-\cdot)-f_{1}(\lambda_{0}-\cdot)\hspace{-0.13cm}\parallel_{p}
\parallel\hspace{-0.13cm}f_{2}\hspace{-0.13cm}\parallel_{p^{'}}\rightarrow 0
\end{equation*}
when $\lambda\rightarrow\lambda_{0}$ for any $\lambda_{0}\in \mathbb{R}^{n}$.
\\
Finally, by taking successive convolutions and using the result of Case 1,
for any $k_{1}>k_{2}>n/\kappa$, $f^{*k_{1}}(\lambda)$ and $f^{*k_{2}}(\lambda)$
are bounded functions, which implies
\begin{equation*}
f^{*k_{1}}(\lambda)\hspace{-0.1cm}=\hspace{-0.1cm}\int_{\mathbb{R}^{n}}\hspace{-0.2cm}
f^{*k_{2}}(\lambda-\eta)f^{*(k_{1}-k_{2})}(\eta)d\eta
\leq \parallel\hspace{-0.1cm} f^{*k_{2}}\hspace{-0.1cm}\parallel_{\infty}
\hspace{-0.15cm}\int_{\mathbb{R}^{n}}\hspace{-0.2cm}f^{*(k_{1}-k_{2})}(\eta)d\eta
=\parallel\hspace{-0.1cm}f^{*k_{2}}\hspace{-0.1cm}\parallel_{\infty}.
\end{equation*}
\qed

%
%
\bigskip
\bibliographystyle{plain}
\begin{small}

\end{small}

\end{document}